\magnification=\magstep1
\input amstex
\documentstyle{amsppt}
\hoffset=.25truein
\hsize=6truein
\vsize=8.95truein

\topmatter
\centerline{\bf A RESIDUE SCALAR PRODUCT FOR ALGEBRAIC}   
\title
FUNCTION FIELDS OVER A NUMBER FIELD
\endtitle
\author
XIAN-JIN LI 
\endauthor
\abstract
    In 1953 Peter Roquette gave an arithmetic proof of the Riemann 
hypothesis for algebraic function fields over a finite constants field, 
which was proved by Andr\'e Weil in 1940.  The construction of 
Weil's scalar product is essential in Roquette's theory.  In this 
paper a scalar product for algebraic function fields over a number 
field is constructed which is the analogue of Weil's scalar product.
\endabstract
\thanks
Research supported by the American Institute of Mathematics.
\endthanks
\address
American Institute of Mathematics, 360 Portage Avenue, Palo Alto, CA 94306
\endaddress
\email
xianjin\@math.Stanford.EDU
\endemail
\endtopmatter
\address
Current address: Department of Mathematics, Brigham Young University, 
Provo, Utah 84602
\endaddress
\document

\heading
1.  Introduction 
\endheading

    In 1953, Roquette [8] gave an arithmetic proof of the 
Riemann hypothesis for algebraic function fields over a finite 
constants field, which was proved by Andr\'e Weil in 1940 [9].  
Since the analogy between number fields and function fields 
over a finite constants field is rather striking,  a question 
is to generalize Roquette's theory to number fields.  
Since the analogue of Roquette's theory 
for number fields does not exist, we want to find the 
analogue of Roquette's theory for function fields over number 
fields.   The construction of Weil's scalar product is essential 
in Roquette's theory.  In this paper a scalar product for 
function fields over number fields is constructed which is the 
analogue of Weil's scalar product.
 
     This paper is organized as follows: In Section 2, 
a set of prime divisors is given on which we base 
our constructions.  In Section 3, we establish some properties 
about the set of multiples of a divisor (cf. [5, Chapter 24]).
These properties are essential for our definition of divisor residues.   
Divisor residues are then defined in Section 4, and they are closely 
related to the theory of correspondences.  Divisor residues may not 
lie in the given  function fields.  In Section 5, we use the 
norm to obtain divisors 
of the given function field from divisor residues.  Next in Section 6, 
we consider function fields over complex numbers which have the 
structure of compact Riemann surfaces.  In this section we make use 
of the Arakelov theory to complete our definition of a residue scalar 
product for algebraic function fields over a number field. 

   Assume that $k$ is a number field.  By the field of algebraic
functions of one variable over $k$ we mean a field, which contains a
transcendental element $x$ and is a finite algebraic extension of the
rational function field $k(x)$.  Let $K$ and $K^\prime$ be two 
function fields of one variable which are finite extensions of
$k(x)$ and $k(x^\prime)$, respectively, with $x$ and $x^\prime$ 
being two algebraically independent elements over $k$.  
The double field $\Delta$ of $K$ and $K^\prime$ is defined to 
be the field $k(x, x^\prime; u, u^\prime)$ with
$f(x, u)=0$ and $f^\prime(x^\prime, u^\prime)=0$ over $k$, where 
$f$ and $f^\prime$ are generating irreducible polynomials of $K$ and 
$K^\prime$ over $k$, respectively.    
Valuations of $k(x, x^\prime)$, which are 
derived from irreducible polynomials of 
$k[x, x^\prime]$, from the negative degree 
of $x$, and from the negative degree of $x^\prime$, correspond
to prime divisors of the rational function field $k(x, x^\prime)$.
Our construction of a residue scalar product is based on those 
prime divisors of $\Delta$ lying over the above prime 
divisors of $k(x, x^\prime)$.  A characterization for these prime 
divisors of the double field $\Delta$ is given in Theorem 2.1.

   To proceed our construction, we start with the multiple ideal of
a divisor, which is a extremely useful tool for this paper.   
Let $\frak o$ be a prime divisor of $K$, and let $\frak o^\prime$ 
be a prime divisor of $K^\prime$.   If $\frak A$ is a divisor 
of the double field $\Delta$, 
the multiple ideal $[\frak A]_{(\frak o, \frak o^\prime)}$ 
is the set of all
elements $a$ of $\Delta$ such that 
$w_{\frak o}(a)\geqslant w_{\frak o}(\frak A)$,
$w_{\frak o^\prime}(a)\geqslant w_{\frak o^\prime}(\frak A)$, 
and $w_{\frak m}(a) \geqslant w_{\frak m}(\frak A)$ 
for all prime divisors $\frak m$ of $\Delta$ lying 
over those prime divisors of $k(x, x^\prime)$ which 
correspond to two variable 
irreducible polynomials of $k[x, x^\prime]$.  
Assume that $\tilde{K^\prime}$
is a finite extension of $K^\prime$.   
Let $\tilde{\frak A}$ be the extension 
of $\frak A$ in the double field of $K, \tilde{K^\prime}$, and let 
$\tilde\frak o^\prime$ be a prime divisor of $\tilde{K^\prime}$ lying
over $\frak o^\prime$.  A relation between the multiple ideals 
$[\frak A]_{(\frak o, \frak o^\prime)}$ and 
$[\tilde{\frak A}]_{(\frak o, \tilde\frak o^\prime)}$ 
is obtained in Theorem 3.3, which will be used in later proofs.

   Let $\frak n$ be a prime divisor of the double field $\Delta$, and
let $\frak A$ be a divisor of $\Delta$, prime to $\frak n$.  
If $\bar{\frak o}$ be a prime divisor of the residue field 
$\Delta\frak n$ of 
$\Delta$ modulo $\frak n$, we denote 
$w_{\bar{\frak o}}(\frak A\frak n)=
\text{min}_{a\in [\frak A]_{\underline{\mu}\bar{\frak o}}}
w_{\bar{\frak o}}(a\frak n)$, where $\underline{\mu}\bar{\frak o}$ is a
pair of prime divisors of $K, K^\prime$ obtained from $\bar{\frak o}$.
A precise formula for $\underline{\mu}\bar{\frak o}$ is 
given in Section 4.   
Then the divisor residue of $\frak A$ modulo $\frak n$ is defined by
$$\frak A\frak n=\sum w_{\bar{\frak o}}(\frak A\frak n)\bar{\frak o} $$
where the sum is over all prime  divisors $\bar{\frak o}$
of $\Delta\frak n$.  It is easy to see that if $\frak A$ is principal,
so is $\frak A\frak n$.  If $\frak A_1$ and $\frak A_2$ are two divisors
of $\Delta$, prime to $\frak n$, we prove in Theorem 4.8 that
$(\frak A_1+\frak A_2)\frak n=\frak A_1\frak n+\frak A_2\frak n$.
The divisor residues are closely related to 
the theory of correspondence,
and a relation between them is given in Theorem 4.11.  

 Let $\frak A$ and $\frak b$ be divisors of $\Delta$, relatively
prime to each other.  Write $\frak A=\sum a_{\frak m}\frak m$ and
$\frak b=\sum b_{\frak n}\frak n$ as linear combination of prime
divisors of $\Delta$ with integer coefficients.  Define
$$\langle\frak A,\frak b\rangle_{f^\prime}=\sum_{\frak m, \frak n}
a_{\frak m}b_{\frak n}\frak N_{\Delta\frak n/K^\prime\nu^\prime}
(\frak m\frak n){\nu^\prime}^{-1}$$
where $\frak N_{L/F}$ is the usual notation of norm and $\nu^\prime$
is the isomorphism of $K^\prime$ induced by $\frak n$.  Then 
$\langle\frak A,\frak b\rangle_{f^\prime}$ is a divisor of $K^\prime$.
By using Theorem 4.8 and Theorem 4.11, we prove in Theorem 5.3
that $\langle\frak A,\frak b\rangle_{f^\prime}
=\langle\frak b, \frak A\rangle_{f^\prime}$.

   When $K^\prime$ is considered as an algebraic function field 
of one variable over complex numbers, 
the set of all places of $K^\prime$
has the structure of a compact Riemann surface
which is denoted by $K^\prime_\infty$.   
If $\frak A$ and $\frak b$ are divisors of $\Delta$, 
we can obtain an Arakelov divisor 
$\langle\frak A,\frak b\rangle_{K^\prime}$ from the divisor
$\langle\frak A,\frak b\rangle_{f^\prime}$ by using the extension
of $\langle\frak A,\frak b\rangle_{f^\prime}$ to $K^\prime_\infty$; 
see Section 6.   Let $\frak d_{K^\prime}$ be a canonical Arakelov
divisor of $K^\prime$.  Define
$$\langle\frak A, \frak b\rangle
=\left(\langle\frak A,\frak b\rangle_{K^\prime}
\cdot \frak d_{K^\prime}\right)$$
where $(\,\cdot\,)$ is the Arakelov intersection product.
Let $\frak A|K^\prime$ be the Arakelov divisor of $K^\prime$ 
obtained from the restriction of $\frak A$ to $K^\prime$.  If
$\frak m$ and $\frak n$ are prime divisors of $\Delta$, we define
$$\{ \frak m, \frak n\}= N_{\Delta/K^\prime}(\frak n)
\left(\frak m|K^\prime\cdot \frak d_{K^\prime}\right)
+N_{\Delta/K^\prime}
(\frak m)\left(\frak n|K^\prime\cdot \frak d_{K^\prime}\right)$$
where $N_{\Delta/K^\prime}(\frak n)=[\Delta\frak n:
K^\prime\nu^\prime]$.  For general divisors $\frak A, \frak b$
of $\Delta$, we define $\{ \frak A, \frak b\}$ by linearity.
Then a residue scalar product of $\frak A, \frak b$ is 
defined by
$$\langle \frak A, \frak b\rangle_r
=\{ \frak A, \frak b\}-\langle \frak A, \frak b\rangle$$
for all divisors $\frak A, \frak b$ of $\Delta$.
By using results of Section 2--5, we proved in Section 6
that the residue scalar product is bilinear and symmetric.  
It is well-defined on the classes of divisors of 
$\Delta$ modulo principal divisors.   

 The author wishes to thank Brian Conrey for his encouragement 
during preparation of the manuscript.  

\heading
2.  Prime divisors
\endheading

     Let $\mu$, $\mu^\prime$ be a pair of isomorphisms 
of $K$, $K^\prime$ into
an algebraic function field $\tilde{K}$ of one variable.  
Then a dependent composite $M$ of $K$, $K^\prime$ is 
defined to be the field composite
$K\mu\cdot K^\prime\mu^\prime$ in $\tilde{K}$.  Consider $M$
as a finite extension of $K^\prime\mu^\prime$.  
Then there exists an isomorphism
$\sigma_0$ of $M$, whose restriction to $K^\prime\mu^\prime$ is
${\mu^\prime}^{-1}$.  Denote $\mu_0=\mu\sigma_0$.  
Then the representative
$M_0=K\mu_0\cdot K^\prime$ of the dependent composite $M$
is called $K^\prime$-normalized.  
If $\sigma$ is an isomorphism of $M_0$ over
$K^\prime$, then $\mu_0\sigma$ is called a 
$K^\prime$-conjugate of $\mu_0$.  
Denote by $\underline{\mu}=\mu_0\underline{\sigma}$ the set of all 
$K^\prime$-conjugates 
of $\mu_0$, and $\underline{\mu}$ is called the isomorphism system 
coordinated to the isomorphism $\mu_0$ of $K$.

   If one of the two isomorphisms $\mu$, $\mu^\prime$, say $\mu$, 
degenerates
to a homomorphism, then a prime ideal $\frak m$ of $K$ exists such that
$\mu$ is the residue class homomorphism modulo $\frak m$ of $K$.
That is, $a\mu$ is the residue $a \frak m$ of $a$ modulo $\frak m$ for
every element $a$ of $K$, where $a\frak m$ represents the symbol $\infty$
when $a$ is not integral for $\frak m$.  The residue field of $K$ modulo 
$\frak m$ is a finite extension of $k$, and the $K^\prime$-normalized 
representative of the dependent composite $M$ of $K\mu$, 
$K^\prime\mu^\prime$ is 
then a finite constants extension of $K^\prime$.

  A binary prime  divisor $\frak m$ of $\Delta$ is defined to be
a non-equivalent normalized valuation $w_{\frak m}$ of $\Delta$
which valuates both $K$ and $K^\prime$ identically zero,  
a $K$-unary prime  divisor $\frak m$ of $\Delta$ is a non-equivalent 
normalized valuation $w_{\frak m}$ of $\Delta$
which valuates $K^\prime$ identically zero, and
 a prime $K^\prime$-unary divisor $\frak m$ of $\Delta$ is a
non-equivalent normalized valuation $w_{\frak m}$ of $\Delta$
which valuates $K$ identically zero.  Denote by $S$ the set of all
binary and unary prime  divisors of $\Delta$.

  In particular, binary prime divisors of the rational function field
$k(x, x^\prime)$ correspond to valuations derived from irreducible 
``binary'' polynomials of $k[x, x^\prime]$, $k(x)$-unary prime 
divisors of $k(x, x^\prime)$ correspond to valuations derived from 
irreducible polynomials of $k[x]$ and to the valuation given by the
negative degree in $x$, and $k(x^\prime)$-unary prime 
divisors of $k(x, x^\prime)$ correspond to valuations derived from 
irreducible polynomials of $k[x^\prime]$ and to the valuation given
by the negative degree in $x^\prime$.  Since $\Delta$ is a finite 
separable extension of $k(x, x^\prime)$, 
all prime divisors of $\Delta$ in $S$ 
are obtained from the corresponding three kinds of prime divisors of
$k(x, x^\prime)$ by extension of valuations under finite separable
extension of fields (cf. [3, \S10]).
The following Lemma gives a characterization of prime
divisors of $\Delta$ in $S$.

\proclaim{Theorem 2.1}  The prime  divisors $\frak m$ of $\Delta$
in $S$ are in one-one correspondence with the classes of 
isomorphism pairs $\mu$, $\mu^\prime$ of $K$, $K^\prime$ into an 
algebraic function
field of one variable in such a way that the residue class homomorphism
of $\Delta$ modulo $\frak m$ induces the isomorphism pair
$\mu$, $\mu^\prime$ and that
$$\Delta\frak m=K\mu\cdot K^\prime\mu^\prime.$$ \endproclaim

  \demo{Proof}  Let $\frak m$ be a prime  divisor of $\Delta$.
Then the residue class homomorphism of $\Delta$ modulo $\frak m$
induces isomorphisms or homomorphisms $\mu$, 
$\mu^\prime$ of $K$, $K^\prime$
into the residue class field $\Delta\frak m$, which is an 
algebraic function field of one variable.  It is clear that 
$\Delta\frak m=K\mu\cdot K^\prime\mu^\prime$, 
and hence, $\frak m$ induces
the class of isomorphism pair $\mu$, $\mu^\prime$ of $K$, $K^\prime$  
satisfying $\Delta\frak m=K\mu\cdot K^\prime\mu^\prime$.

  Conversely, let $\mu$, $\mu^\prime$ be an isomorphism pair of
$K$, $K^\prime$ into an algebraic function field, 
which do not both degenerate
to homomorphisms.  Assume that $\mu^\prime$ does not degenerate to a
homomorphism.  Then the $K^\prime$-normalized representative of the 
dependent composite $K\mu\cdot K^\prime\mu^\prime$ is of the form
$K\mu_0\cdot K^\prime$.  A prime  divisor $\frak m$ of $\Delta$
exists such that $\Delta\frak m=K\mu_0\cdot K^\prime$.  It
follows that the residue homomorphism of $\Delta$
modulo $\frak m$ induces the isomorphism pair $\mu_0$, $1$ of $K$,
$K^\prime$, and $\frak m$ is uniquely determined by the class of 
$K^\prime$-conjugates of $\mu_0$.  Therefore, a unique prime  divisor
$\frak m$ of $\Delta$ exists, which corresponds to the class of
isomorphism pair $\mu$, $\mu^\prime$ of $K$, $K^\prime$ in such 
a way that the residue class homomorphism of $\Delta$ 
modulo $\frak m$ induces $\mu$, $\mu^\prime$ and 
$\Delta\frak m=K\mu\cdot K^\prime\mu^\prime$.

  This completes the proof of the theorem.  \enddemo

  \proclaim{Corollary 2.2}  The $K$-unary and $K^\prime$-unary 
prime  divisors
of $\Delta$ are in one-one correspondence with the prime
 divisors of $K$ and $K^\prime$, respectively.  \endproclaim

  Divisors $\frak A$ of $\Delta$ are defined as the formal sum
$$\frak A=\sum w_{\frak m}(\frak A)\frak m$$
where the sum is over all prime  divisors $\frak m$ in 
$S$ with the integer coefficients $w_{\frak m}(\frak A)$ 
being zero for almost all prime divisors $\frak m$.  

\heading
3.  The ideal of multiples of a divisor
\endheading

 In this section, we study some properties of multiple ideals, 
which are essential for our definition of
divisor residues given in the next section.

  Let $J$ be the set of all elements in $\Delta$ which are integral 
for all binary prime divisors of $\Delta$.  
Let $\frak o$ and $\frak o^\prime$ be
$K$-unary and $K^\prime$-unary prime  divisors of $\Delta$, respectively.
Denote $\underline{\frak o}=(\frak o, \frak o^\prime)$.  The principal
order $J_{\underline{\frak o}}$, restricted by the unary pair 
$\underline{\frak o}$, is the set of all elements of $\Delta$
which are integral for all binary prime  divisors and for the unary pair.

  \proclaim{Lemma  3.1}  $J$ is the ring composite $[K, K^\prime]$ 
of $K$ and $K^\prime$.  \endproclaim

   \demo{Proof}  It is clear that $[K, K^\prime]$ is contained in $J$.
Conversely, let $a$ be a nonzero element of $J$.  Then the denominator
$n$ of $a$ contains only unary prime  divisors.
Let $t$ be a non-constant element of $K$ whose denominator contains
all prime  divisors of $K$ dividing $n$.  Since $K$ is a finite 
algebraic extension of $k(t)$, we can choose a $t$-integral basis
$u_i$ of $K$ over the field $k(t)$.  Then the $u_i$ also
form a basis of $\Delta$ over $K^\prime(t)$.  Since $a$ is integral
in $t$, it can be written in the form
$a=\sum f^\prime_i(t)u_i$ with $f^\prime_i(t)$ in $K^\prime[t]$.  
It follows that $a$ belongs to the ring composite $[K, K^\prime]$. 
\qed \enddemo

  Let $\frak A$ be a  divisor of $\Delta$.  The $J$-multiple ideal 
$[\frak A]$ of $\frak A$ is the set of all elements $a$ in 
$\Delta$ such that $w_{\frak m}(a)\geq w_{\frak m}(\frak A)$
for all binary prime  divisors $\frak m$ of $\Delta$.
  Let $\tilde {K^\prime}$ be a finite extension of $K^\prime$.  Put
$\tilde J=[K, \tilde{K^\prime}]$.  
Let $\tilde{\frak A}$ be the extension of a 
 divisor $\frak A$ of $\Delta$ in the double field $\tilde{\Delta}$
of $K$ and $\tilde{K^\prime}$.   Then the $\tilde J$-multiple ideal
$[\tilde{\frak A}]$ of $\tilde{\frak A}$ in $\tilde{\Delta}$
is $[\frak A]\tilde J$.  In fact, since the extension $\tilde\Delta$ 
of $\Delta$ can be considered as the extension $\tilde {K^\prime}$ 
of $K^\prime$ over $K$, every binary prime divisor
of $\Delta$ splits into  binary prime divisors
of $\tilde\Delta$.  Since binary prime divisors of $\Delta$ form
a proper subset of prime divisors of $\Delta/K$, considered as an
algebraic function field of one variable over $K$, by strong
approximation theorem (cf. [5, Chapter 24])
we find that $\tilde{\frak A}$ is the greatest common binary
divisor of all elements in the $\tilde J$-ideal $[\frak A]\tilde J$.
Since binary prime divisors of $\tilde\Delta$ form
a proper subset of prime divisors of $\tilde\Delta/K$, by strong
approximation theorem we also find that $\tilde{\frak A}$ is the 
greatest common binary divisor of all elements in the $\tilde J$-ideal 
$[\tilde{\frak A}]$.  Therefore, 
we have $[\frak A]\tilde J=[\tilde{\frak A}]$.

  Let $\frak o$, $\frak o^\prime$ be prime  divisors of $K$, $K^\prime$,
and let $1_{\frak o}$, $1^\prime_{\frak o^\prime}$ be the sets
of all elements in $K$, $K^\prime$ which are integral for $\frak o$,
$\frak o^\prime$, respectively.

   \proclaim{Lemma  3.2}   $J_{\underline{\frak o}}=
[1_{\frak o}, 1^\prime_{\frak o^\prime}]$. \endproclaim

  \demo{Proof}   It is clear that 
$[1_{\frak o}, 1^\prime_{\frak o^\prime}]$ 
is contained in $J_{\underline{\frak o}}$.  Conversely, let $a$ be an
element in $J_{\underline{\frak o}}$.  
Then the denominator $n$ of $a$ contains only unary 
prime  divisors of $\Delta$ different from $\frak o$, 
$\frak o^\prime$.  Let $t$ be a non-constant $\frak o$-integral element
in $K$  whose denominator contains all prime  divisors of $K$
dividing $n$.  Choose a $t$-integral basis $u_i$ for $K$
over $k(t)$.  Then the $u_i$ are $\frak o$-integral.  Write
$$a=\sum_i\sum_{k\geq 0}a_{ik}^\prime t^k u_i$$
with $a_{ik}^\prime$ in $K^\prime$.  This implies that $a$ belongs to
$[1_{\frak o}, K^\prime]$.

   Since $t^k u_i$ is in $K$, we have $w_{\frak o^\prime}(t^k u_i)=0$.
Hence, $w_{\frak o^\prime}(a)\geq\text{min}_{i, k}\,
w_{\frak o^\prime}(a^\prime_{ik})=\ell$.  We want to show that $\ell$
is nonnegative.   Arguing by contradiction, assume that $\ell$ is a
negative integer.  Let $\pi^\prime$ be a prime element of $\frak o^\prime$ 
in $K^\prime$.  
Define $b^\prime_{ik}={\pi^\prime}^{-\ell}a^\prime_{ik}$.  Since
$a$ is integral for $\frak o^\prime$, we have
$$\sum_i\sum_{k\geq 0}b_{ik}^\prime t^k u_i\equiv 0\,\, \, 
\text{mod}\, \, \frak o^\prime.$$
Replacing the $b^\prime_{ik}$ by its residues modulo $\frak o^\prime$,
we obtain an equation of the $t^ku_i$ with coefficients in the
residue class field of $K^\prime$ modulo $\frak o^\prime$.
Since the $t^k u_i$ are linearly independent over $k$,
by the argument in [2, Chapter 15] they are linearly independent over the 
residue class field of $K^\prime$ modulo $\frak o^\prime$.  
It follows that
these coefficients must be all equal to zero.  This is impossible.  
Therefore, $\ell\geq 0$, and hence $a$ belongs to $[1_{\frak o},
1^\prime_{\frak o^\prime}]$.\qed\enddemo

  Let $\frak A$ be a  divisor of $\Delta$ and 
$\underline{\frak o}=(\frak o, \frak o^\prime)$ a unary pair.  The 
$J_{\underline {\frak o}}$-multiple ideal 
$[\frak A]_{\underline{\frak o}}$
of $\frak A$ is defined to be the set of all elements $a$ in $\Delta$
such that $w_{\frak m}(a)\geq w_{\frak m}(\frak A)$ for all binary prime
 divisors $\frak m$ of $\Delta$ and such that
$w_{\frak o}(a)\geq w_{\frak o}(\frak A)$,
$w_{\frak o^\prime}(a)\geq w_{\frak o^\prime}(\frak A)$.  
 Let $\tilde{K^\prime}$ be a finite extension
of $K^\prime$.  Assume that $\tilde{\frak o^\prime}$ 
is a prime  divisor of
$\tilde{K^\prime}$ lying above $\frak o^\prime$.  Denote 
$\tilde {\underline{\frak o}}=(\frak o, \tilde{\frak o^\prime})$.
Define 
$\tilde {1^\prime}_{\tilde{\frak o^\prime}}$
to be the set of all $\tilde{\frak o^\prime}$-integral elements of
$\tilde {K^\prime}$.  Then 
$\tilde J_{\tilde{\underline{\frak o}}}=[1_{\frak o},
\tilde {1^\prime}_{\tilde{\frak o^\prime}}]$.

 \proclaim{Theorem 3.3}  Let $\frak A$ be a  divisor of $\Delta$.
Then the $\tilde J_{\tilde{\underline{\frak o}}}$-multiple ideal 
$[\tilde{\frak A}]_{\tilde{\underline{\frak o}}}$ 
of $\tilde{\frak A}$ in $\tilde{\Delta}$ is  
$[\frak A]_{\underline{\frak o}}
\tilde J_{\tilde{\underline{\frak o}}}$. 
\endproclaim

  \demo{Proof}   Multiplying $\frak A$ by an element of $\Delta$,  
we can assume that $\frak A$ is integral for all binary prime
 divisors of $\Delta$.  It will first be shown that
$$[\tilde{\frak A}]_{\underline{\frak o}}=
[\frak A]_{\underline{\frak o}} \tilde J_{\underline{\frak o}},$$
where an element $\tilde a$ of $\tilde{\Delta}$ is said be 
integral for $\frak o^\prime$ if $\tilde a\equiv 0$ 
modulo $\frak o^\prime$.
Let  $[\tilde{\frak A}]_{\underline{\frak o}}^0$ be the set of all
$\underline{\frak o}$-multiples of $\tilde{\frak A}$ in 
$\tilde J_{\underline{\frak o}}$.   Then
$[\tilde{\frak A}]_{\underline{\frak o}}=[\tilde{\frak A}]\cap
[\tilde{\frak A}]_{\underline{\frak o}}^0$.  Similarly,
we have
$[\frak A]_{\underline{\frak o}}=[\frak A]\cap
[\frak A]_{\underline{\frak o}}^0$.  Let $\{\tilde{u_i^\prime}\}$ be an 
$\frak o^\prime$-integral basis for $\tilde{K^\prime}$ 
over $K^\prime$.  Then
$\tilde{K^\prime}=\sum K^\prime \tilde{u_i^\prime}$ and
$\tilde{1^\prime}_{\frak o^\prime}=\sum 1^\prime_{\frak o^\prime}
\tilde {u_i^\prime}$.  It follows that $\tilde J=
[K, \tilde{K^\prime}]=\sum J\tilde{u_i^\prime}$ and
$\tilde J_{\underline{\frak o}}=[1_{\frak o},
\tilde{1^\prime}_{\frak o^\prime}]=\sum J_{\underline{\frak o}}
\tilde{u_i^\prime}$.  This implies that $[\tilde{\frak A}]
=[\frak A]\tilde J=\sum [\frak A]\tilde{u_i^\prime}$.  Let $a$,
$a^\prime$ be elements of $K$, $K^\prime$ such that $w_{\frak o}(a)
=w_{\frak o}(\frak A)$ and $w_{\frak o^\prime}(a^\prime)=
w_{\frak o^\prime}(\frak A)$.  Then 
$[\tilde{\frak A}]_{\underline{\frak o}}^0
=a a^\prime\tilde J_{\underline{\frak o}}=\sum aa^\prime 
J_{\underline{\frak o}}\tilde{u_i^\prime}=\sum
[\frak A]_{\underline{\frak o}}^0 \tilde{u_i^\prime}$.
It follows that
$[\tilde{\frak A}]_{\underline{\frak o}}=[\tilde{\frak A}]\cap
[\tilde{\frak A}]_{\underline{\frak o}}^0=\sum ([\frak A]\cap
[\frak A]_{\underline{\frak o}}^0)\tilde{u_i^\prime}=
\sum [\frak A]_{\underline{\frak o}}\tilde{u_i^\prime}=
[\frak A]_{\underline{\frak o}}\tilde J_{\underline{\frak o}}$.
Therefore, in order to prove the stated result it suffices to show that
$[\tilde{\frak A}]_{\tilde{\underline{\frak o}}}
=[\tilde{\frak A}]_{\underline{\frak o}}
\tilde{1^\prime}_{\tilde{\frak o^\prime}}$. 
 It is clear that $[\tilde{\frak A}]_{\underline{\frak o}}
\tilde{1^\prime}_{\tilde{\frak o^\prime}}$ is contained in
$[\tilde{\frak A}]_{\tilde{\underline{\frak o}}}$.  Conversely,
let $\tilde a$ be an element of $\tilde{\Delta}$ in
$[\tilde{\frak A}]_{\tilde{\underline{\frak o}}}$.  
By strong approximation theorem, we can choose an element 
$\tilde {n^\prime}$ in $\tilde{K^\prime}$ such that
$w_{\tilde{\frak o^\prime}}(\tilde{n^\prime})=0$ and
$w_{\tilde{\frak o^\prime_i}}(\tilde{n^\prime})\geq
w_{\tilde{\frak o^\prime_i}}(\frak A)
-w_{\tilde{\frak o^\prime_i}}(\tilde a)$
for every prime  divisor $\tilde{\frak o^\prime_i}\neq \tilde{\frak
o^\prime}$ of $\tilde{K^\prime}$ lying above $\frak o^\prime$.
Then $a\tilde{n^\prime}$ belongs to
$[\tilde{\frak A}]_{\underline{\frak o}}$, and
$\tilde{n^\prime}^{-1}$ belongs to 
$\tilde {1^\prime}_{\tilde{\frak o^\prime}}$.  Therefore,
$\tilde a$ belongs to 
$[\tilde{\frak A}]_{\underline{\frak o}}
\tilde{1^\prime}_{\tilde{\frak o^\prime}}$.

  This completes the proof of the theorem.  \enddemo

   Let $\frak m$ be a prime  divisor of $\Delta$ 
with $\mu$, $\mu^\prime$ being the isomorphism pair of $K$, $K^\prime$
induced by $\frak m$.  If $\frak m$ is not $K^\prime$-unary, then
degree of $\frak m$ over $K^\prime$ is defined by
$N_{\Delta/K^\prime}(\frak m)=[\Delta\frak m:
K^\prime\mu^\prime]$.  And if $\frak m$ is not $K$-unary, then
the degree of $\frak m$ over $K$ is defined by
$N_{\Delta/K}(\frak m)=[\Delta\frak m:
K \mu]$.   Define $N_{\Delta/K^\prime}(\frak m)=0$
if $\frak m$ is $K^\prime$-unary, and define $N_{\Delta/K}(\frak
m)=0$ if $\frak m$ is $K$-unary.
By a splitting field for a binary prime  divisor $\frak m$ we mean a finite
extension $\tilde{K^\prime}$ of $K^\prime$ such that, in the double
field $\tilde{\Delta}$ of $K$ and $\tilde{K^\prime}$,
$\frak m$ splits into prime  divisors $\tilde{\frak m}$ of degree one
over $\tilde{K^\prime}$.

   Let $\frak m$ be a binary prime  divisor of $\Delta$. 
Assume that the residue class field $\Delta\frak m$ is
$K^\prime$-normalized and that $\underline{\mu}=\mu\underline{\sigma}$
is the coordinated isomorphism system.  Let $\tilde{K^\prime}$ be
a finite extension of $K^\prime$.  Decompose $\mu\underline{\sigma}$
into $\tilde{K^\prime}$-conjugate classes with the $\sigma_i$ being 
representatives for these classes.  By Theorem 2.1, let
$\tilde{\frak m_i}$ be the prime  divisor of $\tilde{\Delta}$
corresponding to the $\tilde{K^\prime}$-conjugate class represented by
$\mu\sigma_i$.  Then the degree of $\tilde{\frak m_i}$ over
$\tilde{K^\prime}$ is
$N_{\tilde{\Delta}/\tilde{K^\prime}}(\tilde{\frak m_i})=
[K\mu\sigma_i\cdot\tilde{K^\prime}:\tilde{K^\prime}]$.  It follows that
a finite extension $\tilde{K^\prime}$ of $K^\prime$ is a splitting field of
$\frak m$ if and only if $K\mu\sigma$ is contained in $\tilde{K^\prime}$
for every isomorphism $\sigma$ of $\Delta\frak m$, and that
the number of such different prime  divisors $\tilde{\frak m_i}$
of $\tilde{\Delta}$ lying above $\frak m$ is equal to the field extension
degree of $\Delta\frak m$ over $K^\prime$.

 \proclaim{Lemma 3.4}  Let $\tilde{K^\prime}$ be a splitting field for 
a binary prime  divisor $\frak m$ of $\Delta$.  Then
$\frak m$ is unramified in the double field $\tilde{\Delta}$
of $K$ and $\tilde{K^\prime}$.  \endproclaim

  \demo{Proof}   Let $\underline{u^\prime}=\{\tilde{u_i^\prime}\}$
be a basis of $\tilde{K^\prime}$ over $K^\prime$.  
Then $\underline{u^\prime}$ is also a $\frak m$-integral basis
for $\tilde{\Delta}$ over $\Delta$.  In fact,
let $\tilde a$ be an $\frak m$-integral element in 
$\tilde{\Delta}$.  Multiply $\tilde a$ by an element 
$\tilde v$ of $\tilde{\Delta}$, which is prime to
$\frak m$, so that $\tilde a\tilde v$ is integral for every binary
prime  divisor of $\tilde{\Delta}$.  Then $\tilde a\tilde v$
belongs to $\tilde J$.  Since $\tilde J=[K, \tilde{K^\prime}]=
\sum J \tilde{u_i^\prime}$,  $\tilde a$ belongs
to $\sum \tilde v^{-1} J\tilde {u_i^\prime}$, and hence
$\tilde a$ can be expressed as a linear combination of $\tilde{u_i^\prime}$
with $\frak m$-integral coefficients.  Since $\tilde {K^\prime}$ is a
splitting field of $\frak m$, there exists no $\frak m$-integral
elements $a_i$ of $\Delta$ such that $\sum a_i \tilde {u_i^\prime}=0$.
Therefore, $\underline{u^\prime}$ is a $\frak m$-integral basis
for $\tilde{\Delta}$ over $\Delta$.  

   Since $\tilde{K^\prime}$ is separable
over $K^\prime$ and since $\underline{u^\prime}$ is a basis 
for $\tilde {K^\prime}$ over $K^\prime$, the discriminant 
$d_{\tilde{K^\prime}/K^\prime}(\underline{u^\prime})$
does not vanish.  
Since $d_{\tilde{K^\prime}/K^\prime}(\underline{u^\prime})$
is an element of $K^\prime$, $\frak m$ is not a divisor of it.
 Now, the discriminant $d_{\tilde{\Delta}/\Delta}(\underline{u^\prime})$
divides  $d_{\tilde{K^\prime}/K^\prime}(\underline{u^\prime})$, and 
hence $\frak m$ is not a  divisor of
$d_{\tilde{\Delta}/\Delta}(\underline{u^\prime})$.
On the other hand, the contribution of $\frak m$ to field discriminant
$d_{\tilde{\Delta}/\Delta}$ of  $\tilde{\Delta}$ over $\Delta$ divides
$d_{\tilde{\Delta}/\Delta}(\underline{u^\prime})$.
Therefore, $\frak m$ is not a  divisor of
$d_{\tilde{\Delta}/\Delta}$.  This implies
that $\frak m$ is unramified in $\tilde{\Delta}$
by the Dedekind discriminant theorem. \qed \enddemo

The following result, which follows immediately from Lemma 3.4
and from the argument preceding Lemma 3.4, is a convenient technique
which will be used in later proofs.

  \proclaim{Corollary 3.5}  Let $\frak m$ be a binary prime  divisor 
with $N_{\Delta/K^\prime}(\frak m)=n$, and let
$\tilde{K^\prime}$ be a splitting field for $\frak m$.  Assume that
the residue class field $\Delta\frak m$ is $K^\prime$-normalized
and that $\mu\underline{\sigma}$ is the coordinated isomorphism system of 
$K$ into $\tilde{K^\prime}$.  Then $\frak m$ has in $\tilde{\Delta}$
the decomposition
$$\frak m=\sum\tilde{\frak m_{\sigma}}$$
into $n$ prime  divisors $\tilde{\frak m_\sigma}$ with 
$N_{\tilde{\Delta}/\tilde{K^\prime}}(\tilde{\frak m_\sigma})=1$
and with the coordinated isomorphism $\mu\sigma$ of $K$ into 
$\tilde{K^\prime}$.  \endproclaim

\heading
4.  Definition of  divisor residues 
\endheading

  We first introduce a notation.  If a field $L$ is a finite separable
extension of a field $F$, and if $\frak m$ is a prime divisor
of $F$ whose extension in $L$ has a factorization of the form
$\frak P_1^{e_1}\frak P_2^{e_2}\cdots\frak P_r^{e_r}$, then we define
$$\Pi_{L/F}(\frak P_i)=\frak m$$
for $i=1,\cdots, r$, and extend it to the group of all divisors of
$L$ by additivity.  In other words, $\Pi_{L/F}(\frak P)$ is a prime
divisor of $F$ lying below the prime divisor $\frak P$ of $L$.

  Let $\frak n$ be a prime  divisor of $\Delta$ in $S$.
Assume that $\underline{\mu}=(\mu, \mu^\prime)$ is the isomorphism
pair of $K$, $K^\prime$ induced by $\frak n$ as in Theorem 2.1.  For every
prime  divisor $\bar {\frak o}$ of $\Delta\frak n$,
define
$$ \mu\bar{\frak o}=\cases \frak n, &\text{if $\frak n$ is
$K$-unary}; \\ \Pi_{\Delta\frak n/K\mu}(\bar{\frak o})
\mu^{-1}, &\text{if $\frak n$ is not $K$-unary}. \endcases$$
and
$$ \mu^\prime\bar{\frak o}=\cases \frak n, &\text{if $\frak n$ is
$K^\prime$-unary}; \\ \Pi_{\Delta\frak n/K^\prime\mu^\prime}(\bar{\frak o})
(\mu^\prime)^{-1}, &\text{if $\frak n$ is not $K^\prime$-unary}. \endcases$$
Put $\underline{\mu}\bar{\frak o}=(\mu\bar{\frak o},
\mu^\prime \bar{\frak o})$.

  \proclaim{Lemma 4.1}  For every element $a$ in $K$, we have
$a(\mu\bar{\frak o})= (a\mu)\bar{\frak o}$. \endproclaim

   \demo{Proof}   First we consider the case when $\frak n$ is
a  $K$-unary prime  divisor of $\Delta$.  In this case,
the stated identity can be written as $a\frak n=(a\frak n)\bar{\frak o}$.
If $\Delta\frak n$ is $K^\prime$-normalized, then $\Delta\frak n$
is a finite constants extension of $K^\prime$, and the residue class field 
of $K$ modulo $\frak n$ is contained in the 
constants field of $\Delta\frak n$. 
On the other hand, we know that the residue class
field of $\Delta\frak n$ modulo $\bar{\frak o}$ contains this
constants field as a subfield.  That is, modulo $\bar{\frak o}$ is the
identity map on the constants field.    It follows that
$(a\frak n)\bar{\frak o}=a\frak n$.

   When $\frak n$ is not $K$-unary, it is clear that
$(a\mu)\bar{\frak o}=(a\mu)\Pi_{\Delta\frak n/
K\mu}(\bar {\frak o})=(a\mu)((\mu\bar{\frak o})\mu)$
$=a(\mu\bar{\frak o})$.  \qed  \enddemo

  If $\frak A$ is a  divisor of $\Delta$, 
prime to $\frak n$, then elements of 
$[\frak A]_{\underline{\mu}\bar{\frak o}}$
are $\frak n$-integral for every prime  divisor $\bar{\frak o}$ of
$\Delta\frak n$.  Denote $w_{\bar{\frak o}}(\frak A\frak n)=
\text{min}_{a\in [\frak A]_{\underline{\mu}\bar{\frak o}}}
w_{\bar{\frak o}}(a\frak n)$.   
The {\bf divisor residue} $\frak A\frak n$ of
$\frak A$ modulo $\frak n$ is defined by
$$\frak A\frak n=\sum w_{\bar{\frak o}}(\frak A\frak n)\bar{\frak o} 
\tag 4.2 $$
where the sum is over all prime  divisors $\bar{\frak o}$
of $\Delta\frak n$.  In the following two lemmas we 
are going to show that 
the divisor residue $\frak A\frak n$ is well-defined.

  \proclaim{Lemma 4.3}  For every prime  divisor $\bar{\frak o}$
of $\Delta\frak n$, $J_{\underline{\mu}\bar{\frak o}}\frak n$
is contained in the set $\bar 1_{\bar{\frak o}}$ of all 
$\bar{\frak o}$-integral elements in $\Delta\frak n$. 
\endproclaim

  \demo{Proof}  By Lemma  3.2, we have
$J_{\underline{\mu}\bar{\frak o}}=[1_{\mu\bar{\frak o}}, 
1^\prime_{\mu^\prime\bar{\frak o}}]$.  Then
$J_{\underline{\mu}\bar{\frak o}}\frak n=[1_{\mu\bar{\frak o}}\mu, 
1^\prime_{\mu^\prime\bar{\frak o}}\mu^\prime]$.   It follows from
Lemma 4.1 that elements in $1_{\mu\bar{\frak o}}\mu$ and  
$1^\prime_{\mu^\prime\bar{\frak o}}\mu^\prime$ are 
$\bar{\frak o}$-integral.
Therefore,
$J_{\underline{\mu}\bar{\frak o}}\frak n$
is contained in $\bar 1_{\bar {\frak o}}$. \qed\enddemo

  \proclaim{Lemma  4.4}   Let $\frak n$ be a prime 
divisor of $\Delta$, and 
let $\frak A$ be a  divisor of $\Delta$, prime to $\frak n$.  
Then $w_{\bar{\frak o}}(\frak A\frak n)$
is finite for every prime  divisor $\bar{\frak o}$ of 
$\Delta\frak n$, and is zero for almost all prime  divisors
$\bar{\frak o}$ of $\Delta\frak n$.  \endproclaim
 
   \demo{Proof}  Since $\frak A$ is prime to $\frak n$, 
there exist elements
$u$ and $v$ of $\Delta$, prime to $\frak n$, such that 
$uJ_{\underline{\mu}\bar{\frak o}} \subseteq
[\frak A]_{\underline{\mu}\bar{\frak o}}\subseteq
\frac 1v J_{\underline{\mu}\bar{\frak o}}$.  Hence,
$u\frak n J_{\underline{\mu}\bar{\frak o}}\frak n \subseteq
[\frak A]_{\underline{\mu}\bar{\frak o}}\frak n\subseteq
\frac 1{v\frak n} J_{\underline{\mu}\bar{\frak o}}\frak n$.
By Lemma 4.3, we have
$1\in J_{\underline{\mu}\bar{\frak o}}\frak n\subseteq 
\bar 1_{\bar{\frak o}}$.
It follows that 
$u\frak n \in
[\frak A]_{\underline{\mu}\bar{\frak o}}\frak n\subseteq
\frac 1{v\frak n} \bar 1_{\bar{\frak o}}$.
This implies that $w_{\bar{\frak o}}(u\frak n)
\geq w_{\bar{\frak o}}(\frak A\frak n)\geq -w_{\bar{\frak o}}
(v\frak n)$.  Therefore,
$w_{\bar{\frak o}}(\frak A\frak n)$ is finite.

 Next, there exist elements $u$, $v$ of $\Delta$, prime
to $\frak n$, such that $uJ\subseteq [\frak A]\subseteq \frac 1v J$.
If $\bar{\frak o}$ is a prime  divisor of $\Delta\frak n$
such that $\mu\bar{\frak o}$, $\mu^\prime\bar{\frak o}$ do not occur
in $u$, $v$ and $\frak A$, then
$uJ_{\underline{\mu}\bar{\frak o}} \subseteq
[\frak A]_{\underline{\mu}\bar{\frak o}}\subseteq
\frac 1v J_{\underline{\mu}\bar{\frak o}}$, and hence 
$w_{\bar{\frak o}}(u\frak n)\geq w_{\bar{\frak o}}
(\frak A\frak n)\geq -w_{\bar{\frak o}}(v\frak n)$.  Since
there are only finitely many exceptional prime  divisors $\bar{\frak o}$,
$w_{\bar{\frak o}}(\frak A\frak n)$ is zero  for
almost all prime  divisors $\bar{\frak o}$ of 
$\Delta\frak n$.\qed  \enddemo

  It follows from Lemma 4.4 that the definition of  divisor residues is
well-defined.  

  \proclaim{Lemma  4.5}  Let $\tilde{K^\prime}$ be a finite extension
of $K^\prime$, and let $\frak n$ be a prime  divisor of $\Delta$.
If $\tilde{\frak n}$ is a prime  divisor of $\tilde{\Delta}$
lying above $\frak n$, then $\frak A\tilde{\frak n}=\frak A\frak n$
for every  divisor $\frak A$ of $\Delta$ which is prime to
$\frak n$.  \endproclaim

\demo{Proof}   Let $\tilde{\frak o}$ be any prime  divisor of
$\tilde{\Delta}\tilde{\frak n}$.  Then $\bar{\frak o}=
\Pi_{\tilde{\Delta}\tilde{\frak n}/\Delta\frak n}
(\tilde{\frak o})$ is a prime  divisor of $\Delta\frak n$
lying below $\tilde{\frak o}$.  Let $\underline{\mu}=
(\mu, \mu^\prime)$ be the isomorphism pair of $K$, $K^\prime$ induced
by $\frak n$.  It is clear that $a\tilde{\frak n}=a\frak n$
for every element $a$ in $K$.  Thus, 
$\underline{\tilde{\mu}}=(\mu, \tilde{\mu^\prime})$ is the
isomorphism pair of $K$, $\tilde{K^\prime}$ induced by $\tilde{\frak n}$
as in Theorem 2.1.  It follows that $\mu\tilde{\frak o}=
\mu\bar{\frak o}$.  Then, by Theorem 3.3 we have
$[\tilde{\frak A}]_{\tilde{\underline{\mu}}\tilde{\frak o}} 
=[\frak A]_{\underline{\mu}\bar{\frak o}}
\tilde J_{\tilde{\underline{\mu}}\tilde{\frak o}}$.
Lemma 4.3 says that 
$1\in \tilde J_{\tilde{\underline{\mu}}\tilde{\frak o}}\tilde{\frak n}
\subseteq \tilde 1_{\tilde{\frak o}}$.  It follows that
$w_{\tilde{\frak o}}(\frak A\tilde{\frak n})= \text{min}_{\tilde a\in
[\tilde{\frak A}]_{\tilde{\underline{\mu}}\tilde{\frak o}}}
w_{\tilde{\frak o}}(\tilde a\tilde{\frak n})=\text{min}_{a\in 
[\frak A]_{\underline{\mu}\bar{\frak o}}}
w_{\tilde{\frak o}}(a\tilde{\frak n})
=\text{min}_{a\in [\frak A]_{\underline{\mu}\bar{\frak o}}}
w_{\tilde{\frak o}}(a\frak n)$.  Let
$\bar{\frak o}=\sum e_{\tilde{\frak o}}\tilde{\frak o}$, where
$e_{\tilde{\frak o}}$ is the ramification index of $\tilde{\frak o}$
over $\bar{\frak o}$.  Then $w_{\tilde{\frak o}}(\frak A\tilde{\frak n})
=e_{\tilde{\frak o}}w_{\bar{\frak o}}(\frak A\frak n)$.  Therefore,
$\frak A\tilde{\frak n}=\sum w_{\tilde{\frak o}}(\frak A\tilde{\frak n})
\tilde{\frak o}=\sum w_{\bar{\frak o}}(\frak A\frak n)\bar{\frak o}=
\frak A\frak n$. \qed  \enddemo

  \proclaim{Lemma  4.6}  Let $\frak n$ be a prime  divisor with
$\underline{\mu}=(\mu, \mu^\prime)$ being the induced isomorphism pair 
of $K$, $K^\prime$.   If $\frak A$ is a purely unary  divisor of
$\Delta$, prime to $\frak n$, then $\frak A\frak n$ is  $\frak A\mu$
for a $K$-unary  divisor $\frak A$, and is $\frak A\mu^\prime$ for a 
$K^\prime$-unary  divisor $\frak A$.  \endproclaim

  \demo{Proof}   Let $\frak A$ be a $K$-unary  divisor, say.  If
$\bar{\frak o}$ is any prime  divisor of $\Delta\frak n$,
then $[\frak A]_{\underline{\mu}\bar{\frak o}}=
[\frak A]_{\mu\bar{\frak o}}J_{\underline{\mu}\bar{\frak o}}$.
An element $a$ in $K$ exists such that 
$[\frak A]_{\mu\bar{\frak o}}=a1_{\mu\bar{\frak o}}$.
It follows that
$[\frak A]_{\underline{\mu}\bar{\frak o}}=
aJ_{\underline{\mu}\bar{\frak o}}$, and hence
$[\frak A]_{\underline{\mu}\bar{\frak o}}\frak n=
a\mu J_{\underline{\mu}\bar{\frak o}}\frak n$.
By Lemma 4.3, we know that
$1\in J_{\underline{\mu}\bar{\frak o}}\frak n\subseteq 1_{\bar{\frak o}}$.
This implies that $w_{\bar{\frak o}}(\frak A\frak n)=
w_{\bar{\frak o}}(a\mu)$.  On the other hand, it can be seen that
$w_{\bar{\frak o}}(a\mu)=w_{\bar{\frak o}}(\frak A\mu)$. 
 Therefore, we have $\frak A\frak n=\frak A\mu$.
\qed  \enddemo

  Let $\frak m$ be a prime  divisor of $\Delta$, and let
$\mu$, $\mu^\prime$ be the isomorphism pair of $K$, $K^\prime$ induced by
$\frak m$ as in Theorem 2.1.  Then the prime
correspondence, which corresponds to every prime  divisor
$\frak p^\prime$ ($\neq \frak m$) of $K^\prime$ a  divisor 
$\frak m(\frak p^\prime)$ of $K$, is defined by
$$\frak m(\frak p^\prime)=\cases  
\Pi_{\Delta\frak m/K\mu}(\frak p^\prime\mu^\prime)\mu^{-1},
&\text{if $\frak m$ is binary};\\
\frak m, &\text{if $\frak m$ is $K$-unary};\\
0, &\text{if $\frak m$ is $K^\prime$-unary}.\endcases$$

  \proclaim{Lemma  4.7}   Let $\frak m$ be a prime  divisor of
$\Delta$.  When $\frak m$ is not $K^\prime$-unary, assume that
$\Delta\frak m$ is $K^\prime$-normalized and that $(\mu, \mu^\prime)$ with
$\mu^\prime=1$ is the isomorphism pair of $K$, 
$K^\prime$ induced by $\frak m$.
Let $\frak p^\prime$ be a $K^\prime$-unary prime  divisor, let
$\Delta\frak p^\prime$ be $K$-normalized, and let
$\underline \pi=(\pi, \pi^\prime)$ with $\pi=1$ be the isomorphism
pair of $K$, $K^\prime$ induced by $\frak p^\prime$.  If 
$\frak m(\frak p^\prime)$ is considered as a  divisor of
$\Delta\frak p^\prime$, then 
$\frak m\frak p^\prime=\frak m(\frak p^\prime)$.  
\endproclaim

   \demo{Proof}  When $\frak m$ is $K^\prime$-unary, 
$\frak m\frak p^\prime=0=\frak m(\frak p^\prime)$  by definition.
When $\frak m$ is $K$-unary, $\frak m\frak p^\prime=
\frak m$ by Lemma  4.6, and $\frak m(\frak p^\prime)
=\frak m$ by definition.  Therefore,
$\frak m\frak p^\prime=\frak m(\frak p^\prime)$ for unary
prime  divisor $\frak m$.

  For the remaining of the proof, $\frak m$ is assumed to be binary.
Then the stated identity can be written as
$$\frak m\frak p^\prime=\Pi_{\Delta\frak m/K\mu}(\frak p^\prime)
\mu^{-1}=\mu\frak p^\prime.$$
For every prime  divisor $\bar{\frak o}$ of $\Delta\frak p^\prime$,
we have $\underline \pi\bar{\frak o}=(\pi\bar{\frak o}, \frak p^\prime)$.
Since $\frak m$ is binary, $[\frak m]_{\underline\pi\bar{\frak o}}
\subseteq J_{\underline\pi\bar{\frak o}}$.  Furthermore, 
$a\in [\frak m]_{\underline\pi\bar{\frak o}}$ if, and only if,
$a\in J_{\underline\pi\bar{\frak o}}$ and $a\frak m=0$.

  We first consider the case when $\bar{\frak o}$ does not divide
$\frak m(\frak p^\prime)$, where $\frak m(\frak p^\prime)$ is
considered as a  divisor of $\Delta\frak p^\prime$.
Since $[\frak m]_{\underline\pi\bar {\frak o}}\frak p^\prime\subseteq
J_{\underline\pi\bar{\frak o}}\frak p^\prime\subseteq 
\bar 1_{\bar{\frak o}}$
by Lemma 4.3, we have $w_{\bar{\frak o}}(\frak m\frak p^\prime)\geq 0$.
It is clear that an $\pi\bar{\frak o}$-integral element $u$ 
in $K$ exists such that
$u(\pi\bar{\frak o})\neq 0$ and $u\frak m(\frak p^\prime)=0$.
By Lemma 4.1, $(u\mu)\frak p^\prime=u\frak m(\frak p^\prime)=0$.
It follows that $u\mu$ is $\frak p^\prime$-integral, 
and hence $a=u-u\mu$ belongs
to $J_{\underline\pi\bar{\frak o}}$ and $a\frak m=0$.  Thus, $a$ is
an element in $[\frak m]_{\underline\pi\bar{\frak o}}$, and
$(a\frak p^\prime)\bar{\frak o}=(u\pi-(u\mu)\frak p^\prime)\bar{\frak o}
=u(\pi\bar{\frak o})-(u\frak m(\frak p^\prime))\bar{\frak o}=
u(\pi\bar{\frak o})\neq 0$.  It follows that 
$w_{\bar{\frak o}}(\frak m\frak p^\prime)=0$.

  We now consider the case when $\bar{\frak o}$ is a prime  divisor
of $\Delta\frak p^\prime$ dividing $\frak m(\frak p^\prime)$.
Then $\pi\bar{\frak o}=\frak m(\frak p^\prime)$.  By Lemma  3.2
every element $a$ in $J_{\underline\pi\bar{\frak o}}$ can be written as
$a=\sum u_i\cdot u_i^\prime$ with $u_i$ in $1_{\pi\bar{\frak o}}$ and
$u_i^\prime$ in $1^\prime_{\frak p^\prime}$.  Then $a\frak m=
\sum u_i\mu\cdot u_i^\prime$.  Since
$(u_i\mu)\frak p^\prime=u_i\frak m(\frak p^\prime)=u_i(\pi\bar{\frak o})$
by Lemma 4.1, the $u_i\mu$ are $\frak p^\prime$-integral.  
Hence $a\frak m$
belongs to $J_{\underline\pi\bar{\frak o}}$.  It follows that
 $[\frak m]_{\underline\pi\bar{\frak o}}=\{a-a\frak m:
a\in J_{\underline\pi\bar{\frak o}}\}$.  Thus,
$w_{\bar{\frak o}}(\frak m\frak p^\prime)=
\text{min}_{a\in J_{\underline\pi\bar{\frak o}}}
w_{\bar{\frak o}}((a-a\frak m)\frak p^\prime)$.  By Lemma 4.1,
$(a-a\frak m)\frak p^\prime=\sum(u_i-u_i\mu)\frak p^\prime\cdot
u_i^\prime\frak p^\prime=\sum (u_i-u_i(\pi\bar{\frak o}))\cdot
u_i^\prime\pi^\prime$.  Since the $u_i^\prime\pi^\prime$ lie in a
finite extension of $k$, 
$w_{\bar{\frak o}}(\frak m\frak p^\prime)=w_{\bar{\frak o}}
(\frak m (\frak p^\prime))$ when $\bar{\frak o}|\frak m(\frak p^\prime)$.
 Therefore, we have
$$\frak m\frak p^\prime=\sum w_{\bar{\frak o}}(\frak m\frak p^\prime)
\bar{\frak o}=\sum w_{\bar{\frak o}}(\frak m(\frak p^\prime))
\bar{\frak o}=\frak m(\frak p^\prime).\qed$$
  \enddemo

In the next theorem, we are going to show 
that divisor residues modulo $\frak n$ satisfy the distributive law.

 \proclaim{Theorem 4.8 }  Let $\frak n$ be a prime  divisor of
$\Delta$.  If $\frak A$ and $\frak A_1$ are  divisors of 
$\Delta$, prime to $\frak n$, then
$(\frak A+\frak A_1)\frak n=\frak A\frak n+\frak A_1\frak n$.
If $\frak A$ is integral or principal, so is $\frak A\frak n$. \endproclaim

  \demo{Proof}   (1).  To prove that $(\frak A+\frak A_1)\frak n=
\frak A\frak n+\frak A_1\frak n$, it suffices to show that
$[\frak A+\frak A_1]_{\underline{\frak o}}=[\frak A]_{\underline{\frak o}}
\cdot [\frak A_1]_{\underline{\frak o}}$ for every unary pair
$\underline{\frak o}=(\frak o, \frak o^\prime)$.  In fact, 
for every prime  divisor $\bar{\frak o}$ of $\Delta\frak n$,
there exists an element $a$ in $[\frak A]_{\underline \mu\bar{\frak o}}$
and an element $b$ in $[\frak A_1]_{\underline \mu\bar{\frak o}}$
such that $w_{\bar{\frak o}}(a\frak n)=w_{\bar{\frak o}}(\frak A\frak n)$
and $w_{\bar{\frak o}}(b\frak n)=w_{\bar{\frak o}}(\frak A_1\frak n)$.
This implies that $w_{\bar{\frak o}}((\frak A+\frak A_1)\frak n)
\leq w_{\bar{\frak o}}(\frak A\frak n)+w_{\bar{\frak o}}(\frak A_1\frak n)$.
Conversely, there exists an element $c$ in 
$[\frak A+\frak A_1]_{\underline\mu\bar{\frak o}}$ such that
$w_{\bar{\frak o}}((\frak A+\frak A_1)\frak n)=
w_{\bar{\frak o}}(c\frak n)$.  Since 
$[\frak A+\frak A_1]_{\underline\mu\bar{\frak o}}=
[\frak A]_{\underline\mu\bar{\frak o}}\cdot
[\frak A_1]_{\underline\mu\bar{\frak o}}$,
$c=\sum a_ib_i$ with $a_i$ in
$[\frak A]_{\underline\mu\bar{\frak o}}$ and $b_i$ in
$[\frak A_1]_{\underline\mu\bar{\frak o}}$.
It follows that
$w_{\bar{\frak o}}(c\frak n)\geq w_{\bar{\frak o}}(\frak A\frak n)+
w_{\bar{\frak o}}(\frak A_1\frak n)$.  Therefore,
$w_{\bar{\frak o}}((\frak A+\frak A_1)\frak n)
=w_{\bar{\frak o}}(\frak A\frak n)+w_{\bar{\frak o}}(\frak A_1\frak n)$.
That is,
$(\frak A+\frak A_1)\frak n=\frak A\frak n+ \frak A_1\frak n$.

    Now, we are going to prove that
$[\frak A+\frak A_1]_{\underline{\frak o}}=
[\frak A]_{\underline{\frak o}}\cdot [\frak A_1]_{\underline{\frak o}}$
for every unary pair $\underline{\frak o}$.  Multiplying
$\frak A$ and $\frak A_1$ by suitable nonzero elements of
$\Delta$, we can assume that $\frak A$ and $\frak A_1$
are integral for all binary prime divisors of $\Delta$.
First, it is clear that $[\frak A]_{\underline{\frak o}}\cdot 
[\frak A_1]_{\underline{\frak o}}\subseteq 
[\frak A+\frak A_1]_{\underline{\frak o}}$.
Next, by the fundamental theorem of ideal theory of an algebraic 
function field in [5, Chapter 24], we have
$[\frak A+\frak A_1]_{\frak o^\prime}=
[\frak A]_{\frak o^\prime}\cdot [\frak A_1]_{\frak o^\prime}$
when $\Delta$ is considered as a function field over $K$.
Let $\pi$ be a prime element of $\frak o$ in $K$.  Dividing
$\frak A$ and $\frak A_1$ by some powers of $\pi$, we can assume
that $\frak A$ and $\frak A_1$ are prime to $\frak o$.
In order to prove that $[\frak A+\frak A_1]_{\underline{\frak o}}
\subseteq [\frak A]_{\underline{\frak o}}\cdot 
[\frak A_1]_{\underline{\frak o}}$, it suffices to show that
 $a\in [\frak A]_{\underline{\frak o}}\cdot 
[\frak A_1]_{\underline{\frak o}}$ for every element 
$a\in [\frak A+\frak A_1]_{\underline{\frak o}}$.
Since $[\frak A+\frak A_1]_{\underline{\frak o}}
\subseteq [\frak A+\frak A_1]_{\frak o^\prime}=
[\frak A]_{\frak o^\prime}\cdot [\frak A_1]_{\frak o^\prime}$, we can
write $a=\sum a_i^\prime\cdot b_i^\prime$ with
$a_i^\prime\in [\frak A]_{\frak o^\prime}$ and
$b_i^\prime\in [\frak A_1]_{\frak o^\prime}$.  Then 
a smallest nonnegative integer $k$ exists such 
that $\pi^ka=\sum a_i\cdot b_i$ with $a_i\in 
[\frak A]_{\underline{\frak o}}$
and $b_i\in [\frak A_1]_{\underline{\frak o}}$.  We want to show
that $k=0$.  

   Argue by contradiction, assuming that $k\geqslant 1$.   Let 
$\eta$ be an element of $[\frak A]_{\underline{\frak o}}$ such that
$w_{\frak o^\prime}(\eta\frak o)\leqslant w_{\frak o^\prime}(a_i\frak o)$
for each $i$.  Put $\bar\alpha_i=a_i\frak o/\eta\frak o$.  Then 
$\bar\alpha_i$ is $\frak o^\prime$-integral.   Since
$(a_i-\bar\alpha_i\eta)\frak o=0$, we can write
$a_i-\bar\alpha_i\eta=\pi \bar c_i$ 
for some $\frak o$-integral element $\bar c_i$.
Since $\Delta\frak o$ is a finite constants extension of $K^\prime$,
we denote by $d$ the degree of $\Delta\frak o$ over $K^\prime$.
If $\bar\Delta$ is the double field of $K$ and $\Delta\frak o$,
then $\bar\Delta$ is a finite extension of $\Delta$ of degree $d$.
Let $\alpha_i={1\over d}\text{trace}_{\bar\Delta/\Delta}(\bar\alpha_i)$
and $c_i={1\over d}\text{trace}_{\bar\Delta/\Delta}(\bar c_i)$.
Then $\alpha_i$ is an $\frak o^\prime$-integral element of $\Delta$,
$c_i$ is an $\frak o$-integral element of $\Delta$, and we have
$a_i-\alpha_i\eta=\pi c_i$.
It is clear that $c_i\in  [\frak A]_{\underline{\frak o}}$.
Write $\pi^ka=\pi\sum c_i\cdot b_i+\eta\sum \alpha_i\cdot b_i$.
Since $\eta\frak o\neq 0$, this identity implies that
$\sum \alpha_i\cdot b_i=\pi b$ for some $\frak o$-integral element
$b$.  It is clear that $b\in [\frak A_1]_{\underline{\frak o}}$.
Therefore, we have
$\pi^{k-1}a=\sum c_i\cdot b_i+\eta\cdot b$ with
$c_i, \eta\in [\frak A]_{\underline{\frak o}}$ and
$b_i, b\in [\frak A_1]_{\underline{\frak o}}$.   This contradicts
to the minimality of $k$, and hence we must have $k=0$.
Thus, we have proved that
$[\frak A+\frak A_1]_{\underline{\frak o}}=
[\frak A]_{\underline{\frak o}}\cdot [\frak A_1]_{\underline{\frak o}}$.

 (2).  Assume that $\frak A$ is integral.  
If $\frak A_1$ is a  divisor of $\Delta$ dividing $\frak A$,
then for every prime  divisor $\bar{\frak o}$ of $\Delta\frak n$
we have $[\frak A]_{\underline\mu\bar{\frak o}}\subseteq 
[\frak A_1]_{\underline\mu\bar{\frak o}}$, and hence
$[\frak A]_{\underline\mu\bar{\frak o}}\frak n\subseteq 
[\frak A_1]_{\underline\mu\bar{\frak o}}\frak n$.  This implies that
$w_{\bar{\frak o}}(\frak A\frak n)\geq
w_{\bar{\frak o}}(\frak A_1\frak n)$.  That is,
$\frak A_1\frak n$ divides $\frak A\frak n$.  It follows that
$\frak A\frak n$ is integral.

  (3).  When $\frak A=(a)$, by definition we find that
$\frak A\frak n=(a\frak n)$.  Therefore,
$\frak A\frak n$ is principal.

  This completes the proof of the theorem.  \enddemo

  Finally in this section, we study some properties 
of general correspondence,
which is closely related to divisor residues and will be used in
later proofs.

  Let $\frak A=\sum w_{\frak m}(\frak A)\frak m$ be a  divisor of
$\Delta$, and let $\frak p^\prime$ be a prime  divisor of $K^\prime$, 
prime to $\frak A$.  
Then the correspondence $\frak A(\frak p^\prime)$
of $K^\prime$ in $K$ is defined by
$$\frak A(\frak p^\prime)=\sum_{\frak m} w_{\frak m}(\frak A)
\frak m(\frak p^\prime) \tag 4.9$$
where the sums are over all prime  divisors $\frak m$ of $\Delta$.

  \proclaim{Lemma 4.10}  Let $\frak m$ be a prime  
divisor of $\Delta$, and
let $\tilde{K^\prime}$ be a normal splitting extension for $\frak m$.
If $\tilde{\frak p^\prime}$ is a prime  
divisor of $\tilde{K^\prime}$ lying
above a prime  divisor $\frak p^\prime$ ($\neq \frak m$) 
of $K^\prime$, then
$$\frak m(\tilde{\frak p^\prime})=\frak m(\frak p^\prime).$$ 
\endproclaim

  \demo{Proof}  If $\frak m$ is unary, by definition
$\frak m(\tilde{\frak p^\prime})=\frak m(\frak p^\prime)=\frak m$
or $0$.  If $\frak m$ is binary, assume that
$\Delta\frak m=K\mu\cdot K^\prime$.  Let
$\frak m=\sum \tilde {\frak m_\sigma}$ be the decomposition of
$\frak m$ as given in Corollary 3.5.
We have
$\tilde{\Delta}\tilde{\frak m_\sigma}=\tilde{K^\prime}$.
Let $(\Delta\frak m)\sigma=K\mu\sigma\cdot K^\prime$.  Then
the prime  divisor $\bar{\frak p^\prime_\sigma}$ of
$(\Delta\frak m)\sigma$ containing
$\tilde{\frak p^\prime}$ is given by
$$\bar{\frak p^\prime_\sigma}=\Pi_{\tilde{K^\prime}/
(\Delta\frak m)\sigma} (\tilde{\frak p^\prime}).$$
It follows that
$$\aligned \tilde{\frak m_\sigma}(\tilde{\frak p^\prime}) &=
\Pi_{\tilde{K^\prime}/K\mu\sigma}(\tilde{\frak p^\prime}) 
(\mu\sigma)^{-1}
=\Pi_{(\Delta\frak m)\sigma/K\mu\sigma}
(\bar{\frak p^\prime_\sigma}) \sigma^{-1}\mu^{-1} \\
& =\Pi_{\Delta\frak m/K\mu}(\bar{\frak p^\prime_\sigma}
\sigma^{-1}) \mu^{-1}.  \endaligned $$
Since 
$\sum\bar{\frak p^\prime_\sigma}\sigma^{-1}=
\frak p^\prime$,  we have
$\frak m(\tilde{\frak p^\prime})=\sum_{\underline\sigma}
\tilde{\frak m_\sigma}(\tilde{\frak p^\prime})=
\Pi_{\Delta\frak m/K\mu}(\sum_{\underline\sigma}
\bar{\frak p^\prime_\sigma}\sigma^{-1})\mu^{-1}
=\frak m(\frak p^\prime)$.
  \qed\enddemo

   \proclaim{Theorem 4.11}  Let $\frak p^\prime$ be a prime  divisor
of $K^\prime$, and let $\underline\pi=(\pi, \pi^\prime)$ with $\pi=1$ be 
the isomorphism pair of $K$, $K^\prime$ induced by $\frak p^\prime$.  Then
for any  divisor $\frak A$ of $\Delta$, prime to
$\frak p^\prime$, we have $\frak A\frak p^\prime=\frak A(\frak p^\prime)$
when $\frak A(\frak p^\prime)$ is considered as a  divisor of
$\Delta\frak p^\prime$.
\endproclaim

   \demo{Proof}  By Theorem 4.8  and the definition of general
correspondence, we can assume that $\frak A$ is a prime  divisor 
$\frak m$.  Choose $\tilde{K^\prime}$ to be a normal
splitting extension of $K^\prime$ for $\frak m$.  
Let $\tilde{\frak p^\prime}$
be a prime  divisor of $\tilde{K^\prime}$ 
lying above $\frak p^\prime$.  Then
by Lemma  4.5 and Lemma 4.10 we have 
$\frak m\tilde{\frak p^\prime}=\frak m\frak p^\prime$
and
$\frak m(\tilde{\frak p^\prime})=\frak m(\frak p^\prime)$.
Let $\frak m=\sum \tilde{\frak m_\sigma}$ be the decomposition
of $\frak m$ into primes of degree one or zero over $\tilde{K^\prime}$.
Then $\frak m\tilde{\frak p^\prime}=\sum\tilde{\frak m_\sigma}
\tilde{\frak p^\prime}$ and
$\frak m(\tilde{\frak p^\prime})=\sum\tilde{\frak m_\sigma}
(\tilde{\frak p^\prime})$.  It follows from Lemma  4.7 that
$\tilde{\frak m_\sigma}\tilde{\frak p^\prime}=\tilde{\frak m_\sigma}
(\tilde{\frak p^\prime})$, and hence
$\frak m\tilde{\frak p^\prime}=\frak m(\tilde{\frak p^\prime})$.
Consequently, we have 
$\frak m\frak p^\prime=\frak m(\frak p^\prime)$.\qed
  \enddemo

  \proclaim{Corollary 4.12}  Let $\tilde{K^\prime}$ be a finite extension
of $K^\prime$, and let $\tilde{\frak p^\prime}$ be a prime  divisor of
$\tilde{K^\prime}$ lying above a prime  divisor $\frak p^\prime$ of 
$K^\prime$.  Then for any  divisor $\frak A$ of $\Delta$,
prime to $\frak p^\prime$, we have 
$\frak A(\tilde{\frak p^\prime})=\frak A(\frak p^\prime)$.\endproclaim

  \demo{Proof}  Let $\underline\pi=(\pi, \pi^\prime)$ be the
isomorphism pair of $K$, $K^\prime$ induced  by $\frak p^\prime$.
By Theorem 4.11, we have $\pi\frak A(\tilde{\frak p^\prime})=
\frak A\tilde{\frak p^\prime}$ and $\frak A\frak p^\prime=
\pi\frak A(\frak p^\prime)$.  In addition, by Lemma  4.5 we have
$\frak A\tilde{\frak p^\prime}=\frak A\frak p$.  Hence,
$\pi\frak A(\tilde{\frak p^\prime})=\pi\frak A(\frak p^\prime)$.
Since $\pi$ is an isomorphism,
$\frak A(\tilde{\frak p^\prime})=\frak A(\frak p^\prime)$.
\qed\enddemo

  \proclaim{Corollary 4.13}  Let $\frak p^\prime$ be a prime  divisor of
$K^\prime$, and let $\frak A$, $\frak A_1$ and $\frak A_2$ be  divisors of
$\Delta$, prime to $\frak p^\prime$.  Then

(1). $(\frak A_1+\frak A_2)(\frak p^\prime)=
    \frak A_1(\frak p^\prime)+\frak A_2(\frak p^\prime)$.

(2). If $\frak A$ is integral or principal, so is
$\frak A(\frak p^\prime)$.

(3).  $$ \frak A(\frak p^\prime)=\cases \frak A,
    &\text{if $\frak A$ is $K$-unary};\\
      0, &\text{if $\frak A$ is $K^\prime$-unary}. \endcases$$ 
\endproclaim

  \demo{Proof}  We have (1) and (2) by Theorem 4.11 and
Theorem 4.8 , and we have (3) by (1) and the definition of
correspondence. \qed\enddemo

\heading
5. Norm of divisor residues
\endheading

   In this section, we prove a theorem about the norm of divisor 
residues, which is essential for our definition of a residue 
scalar product for algebraic function fields over number fields.

   Assume that $L$ is a finite separable extension of $F$, and let 
$\frak P$ be a prime divisor of $L$.  Define
$$\frak N_{L/F}(\frak P)=\sum_\sigma \frak P\sigma $$
where the sum is over all isomorphisms of $L$ over $F$.
We extend the map $\frak N_{L/F}$ to the group of
all divisors of $L$ by additivity.  

  \proclaim{Lemma  5.1}  Let $\frak m$ and $\frak n$ be different
binary prime  divisors of $\Delta$, and let 
$\underline\pi=(\pi, \pi^\prime)$ be the isomorphism pair of
$K$, $K^\prime$ induced by $\frak n$.  If $\bar{\frak o}$ is a prime
 divisor of $\Delta\frak n$, then
$\bar{\frak o}|\frak m\frak n$ if, and only if, 
$\pi\bar{\frak o}|\frak m(\pi^\prime\bar{\frak o})$.  \endproclaim

  \demo{Proof}   Denote $\frak o=\pi\bar{\frak o}$ and
$\frak o^\prime=\pi^\prime\bar{\frak o}$.  Then $\frak o$ and
$\frak o^\prime$ are unary prime  divisors.  Moreover, $\frak o$ divides
$\frak m(\frak o^\prime)$ if, and only if, $\frak o\mu$ divides
$\Pi_{\Delta\frak m/K\mu}(\frak o^\prime\mu^\prime)$,
where $\underline\mu=(\mu, \mu^\prime)$ is the isomorphism pair
of $K$, $K^\prime$ induced by $\frak m$.  Since $\frak o\mu$ is a
prime  divisor in $K\mu$, $\frak o\mu$ divides
$\Pi_{\Delta\frak m/K\mu}(\frak o^\prime\mu^\prime)$
if, and only if, $\frak o\mu$ and $\frak o^\prime\mu^\prime$
are not coprime in $\Delta\frak m$.

   We first consider the case when $\frak o$ divides 
$\frak m(\frak o^\prime)$.  Let $\bar{\frak o_{\frak m}}$ be a
common prime  divisor of $\frak o\mu$ and $\frak o^\prime\mu^\prime$
in $\Delta\frak m$.  Then $\underline\mu\bar{\frak o_{\frak m}}=
(\frak o, \frak o^\prime)$.  Since $\frak m$ is a binary prime  divisor,
$[\frak m]_{\underline\pi\bar{\frak o}}$ is the set of all elements
$a$ in $J_{\underline\pi\bar{\frak o}}$ such that
$a\frak m=0$.  Now, for every element $a=\sum u_i\cdot u^\prime_i$ in 
$[\frak m]_{\underline\pi\bar{\frak o}}$ with $u_i\in 1_{\pi\bar{\frak o}}$
and $u_i^\prime\in 1^\prime_{\pi^\prime\bar{\frak o}}$, we have
$(a\frak n)\bar{\frak o}=\sum(u_i\pi)\bar{\frak o}\cdot
(u^\prime_i\pi^\prime)\bar{\frak o}=\sum u_i\frak o
\cdot u_i^\prime\frak o^\prime
=\sum u_i(\mu\bar{\frak o_{\frak m}})\cdot u_i^\prime
(\mu^\prime\bar{\frak o_{\frak m}})=(a\frak m)\bar{\frak o_{\frak m}}=0$.
It follows that $w_{\bar{\frak o}}(\frak m\frak n)\geq 1$.  That is,
$\bar{\frak o}$ divides $\frak m\frak n$.

   We now consider the case when $\frak o$ does not divide 
$\frak m(\frak o^\prime)$.  
Then $\frak o\mu$ and $\frak o^\prime\mu^\prime$
are coprime.  It follows that an element $u$ in $K$ exists such that
$u\mu\equiv 1$ modulo $\frak o\mu$ and $u\mu\equiv 0$ modulo
$\frak o^\prime\mu^\prime$.  Assume that 
$x^k+a_1 x^{k-1}+\cdots+a_k$ be the irreducible generating polynomial of
$\Delta\frak m$ over $K^\prime\mu^\prime$.
Define $\underline u=(u^k, \cdots, u, 1)$ and
$\underline{u^\prime}=(1, a_1, \cdots, a_k)(\mu^\prime)^{-1}$.  Then
$\underline u\mu\cdot \underline{u^\prime}\mu^\prime=0$.  Since
 $u\mu\equiv 1$ modulo $\frak o\mu$ and $u\mu\equiv 0$ modulo
$\frak o^\prime\mu^\prime$, we have 
$\underline u \frak o=(1, 1, \cdots, 1)$ and 
$\underline{u^\prime}\frak o^\prime
=(1, 0, \cdots, 0)$.  Let $a=\underline u\cdot\underline{u^\prime}$.
Then $a$ belongs to $[\frak m]_{\underline\pi\bar{\frak o}}$.
Moreover, $(a\frak n)\bar{\frak o}=\underline u\frak o
\cdot\underline{u^\prime}\frak o^\prime=1$.
This implies that $w_{\bar{\frak o}}(\frak m\frak n)=0$, and therefore, 
$\bar{\frak o}$ does not divide $\frak m\frak n$.\qed  \enddemo

 Let $\underline\mu=(\mu_1, \mu_2)$ be a pair of different
isomorphisms of $K$ into an algebraic function field.
Denote $\tilde{K^\prime}=K\mu_1\cdot K\mu_2$ the field composite of
$K\mu_1$ and $K\mu_2$.  For every prime
 divisor $\tilde{\frak o}$ of $\tilde{K^\prime}$, define
$\ell_{\tilde{\frak o}}=\text{min}_{a\in 1_{\underline\mu\tilde{\frak o}}}
w_{\tilde{\frak o}}(a\mu_1-a\mu_2)$, where
$\underline\mu\tilde{\frak o}=(\mu_1\tilde{\frak o},
\mu_2\tilde{\frak o})$ and $1_{\underline\mu\tilde{\frak o}}$ is the
set of all elements in $K$ which are integral for $\mu_1\tilde{\frak
o}$ and $\mu_2\tilde{\frak o}$.  
A different  divisor $\frak D_{\underline\mu}$
for the isomorphism pair $\underline\mu$ is defined by
$$\frak D_{\underline\mu}=\sum_{\tilde{\frak o}}\ell_{\tilde{\frak o}}
\tilde{\frak o}$$
where the sum is over all prime  divisors $\tilde{\frak o}$
of $\tilde{K^\prime}$.
  In general, if $\underline\mu=(\mu_1, \cdots, \mu_n)$ is a system
of different isomorphisms of $K$ into an algebraic function field,
then the conjugate different  divisors $\frak D^{(i)}_{\underline\mu}$
are defined by
$$\frak D^{(i)}_{\underline\mu}=\sum_{j=1, j\neq i}^m
\frak D_{\mu_i, \mu_j}$$
for $i=1,2,\cdots,n$, and the discriminant  
divisor $\frak D_{\underline\mu}$ 
is then defined by
$$\frak D_{\underline\mu}=\sum_{i=1}^m \frak D^{(i)}_{\underline\mu}.$$

  \proclaim{Lemma 5.2}   Let $\frak m$ and $\frak n$ be two 
different prime  divisors of $\Delta$, which are not $K^\prime$-unary.
Assume that $\Delta\frak m=K\mu_0\cdot K^\prime$
and $\Delta\frak n=K \nu\cdot K^\prime$.  If $\underline\mu$
is the isomorphism system of $K$ coordinated to $\mu_0$, then
$$\frak m\frak n=\sum_{\mu\in\underline\mu}\frak D_{\mu, \nu}$$
where the sum is over all isomorphisms in the
isomorphism system $\underline\mu$.  \endproclaim

   \demo{Proof}  Assume that $\tilde{K^\prime}$ is a splitting extension of
$K^\prime$ for $\frak m$ which contains $K\nu$.  Since $\frak m$
and $\frak n$ are different prime  divisors of $\Delta$,
by Theorem 2.1 $\mu$ is not equal to $\nu$ for every isomorphism
$\mu$ in $\underline\mu$.

  (1).  We first consider the case when $\frak m$ and $\frak n$
are binary prime  divisors of $\Delta$ with
$N_{\Delta/K^\prime}(\frak m)=1$.   In this case, $\underline\mu$
consists of a single isomorphism $\mu$, and $K\mu\subseteq
\Delta\frak m=K^\prime\subseteq \Delta\frak n=
K\nu\cdot K^\prime$.  Hence, we can choose 
$\tilde{K^\prime}=\Delta\frak n$
by the argument preceding the statement of Lemma 3.4.  The
stated identity can then be written as
$$w_{\bar{\frak o}}(\frak m\frak n)=w_{\bar{\frak o}}
(\frak D_{\mu, \nu})$$
for every prime  divisor $\bar{\frak o}$ of $\Delta\frak n$.

  Assume that $\mu\bar{\frak o}$ is not equal to $\nu\bar{\frak o}$.
Then $w_{\bar{\frak o}}(\frak D_{\mu, \nu})=0$ by definition.
In order to see that $w_{\bar{\frak o}}(\frak m\frak n)=0$,
it suffices to show that $\nu\bar{\frak o}$ does not divide
$\frak m(\nu^\prime\bar{\frak o})$ by Lemma  5.1, where $\nu$,
$\nu^\prime$ with $\nu^\prime=1$ is the isomorphism pair of
$K$, $K^\prime$ induced by $\frak n$.   We have
$\frak m(\nu^\prime\bar{\frak o})=\Pi_{\Delta\frak m/
K\mu}(\nu^\prime\bar{\frak o})\mu^{-1}$.
Since $\nu^\prime\bar{\frak o}=\Pi_{\Delta\frak n/K^\prime}
(\bar{\frak o})$, we have
$\frak m(\nu^\prime\bar{\frak o})=\Pi_{\Delta\frak n/
K\mu}(\bar{\frak o})\mu^{-1}=\mu\bar{\frak o}$.
Since $\mu\bar{\frak o}\neq \nu\bar{\frak o}$ by assumption, we
have $\nu\bar{\frak o}$ does not divide 
$\frak m(\nu^\prime\bar{\frak o})$.

  Next, assume that $\mu\bar{\frak o}=\nu\bar{\frak o}$.
Denote $\frak o=\mu\bar{\frak o}=\nu\bar{\frak o}$.  If $u$
is integral for $\frak o$, then $u\mu$ is integral for
$\frak o\mu=\Pi_{\Delta\frak n/K\mu}(\bar{\frak o})
=\Pi_{K^\prime/K\mu}(\nu^\prime\bar{\frak o})$.
This implies that $u\mu$ is integral for $\nu^\prime\bar{\frak o}$.
Hence, elements of the form $u-u\mu$ with $u$ in $1_{\frak o}$
belongs to $[\frak m]_{\underline{\frak o}}$ by Lemma  3.2, where
$\underline{\frak o}=(\frak o, \nu^\prime\bar{\frak o})$.
Since $\frak m$ is a binary prime  divisor,
$[\frak m]_{\underline{\frak o}}$ is the set of elements
$a$ in $J_{\underline{\frak o}}$ such that $a\frak m=0$.
If $a$ belongs to $[\frak m]_{\underline{\frak o}}$,
then $a=\underline u\cdot\underline{u^\prime}$ and $a\frak m=0$
with $\underline u$, $\underline {u^\prime}$ being vectors
of equal length over $1_{\frak o}$, $1_{\nu^\prime\bar{\frak o}}$,
and hence $a=\underline u\cdot\underline{u^\prime}-\underline u\mu\cdot
\underline{u^\prime}=(\underline u-\underline u\mu)\cdot 
\underline {u^\prime}$.  It follows that
$[\frak m]_{\underline{\frak o}}$ is generated over 
$1^\prime_{\nu^\prime\bar{\frak o}}$ by elements of the form
$u-u\mu$ with $u$ in $1_{\frak o}$.  This implies that
$$w_{\bar{\frak o}}(\frak m\frak n)=
\text{min}_{u\in 1_{\frak o}}(u\mu-u\nu)$$
by the definition of  divisor residues.  Therefore, by the
definition of different  divisors we have
$$w_{\bar{\frak o}}(\frak m\frak n)
=w_{\bar{\frak o}}(\frak D_{\mu, \nu}).$$
Thus, we have verified the stated identity in the case when $\frak m$ 
and $\frak n$ are binary prime  divisors of $\Delta$
with $N_{\Delta/K^\prime}(\frak m)=1$.   

(2).   We now consider the general case when $\frak m$
and $\frak n$ are binary prime  divisors of $\Delta$
with $N_{\Delta/K^\prime}(\frak m)$ arbitrary.   
Let $\tilde{K^\prime}$ be a splitting extension for $\frak m$,
 which contains $K\nu$.  Then, by Corollary 3.5, $\frak m$ has
in the double field $\tilde\Delta$ of $K$ and 
$\tilde{K^\prime}$ the decomposition
$$\frak m=\sum_{\underline\mu}\tilde {\frak m}_\mu$$
into prime  divisors $\tilde{\frak m}_\mu$ with
$N_{\tilde\Delta/\tilde{K^\prime}}(\tilde{\frak m}_\mu)=1$.
Let $\tilde{\frak n}$ be a prime  divisor of
$\tilde \Delta$ lying above $\frak n$.  By Lemma  4.5, 
we have
$\frak m\frak n=\frak m\tilde{\frak n}$.
By using Theorem 4.8, we obtain that
$$\frak m\frak n=\sum_{\underline\mu}\tilde{\frak m}_\mu
\tilde{\frak n}.$$
Note that $\tilde{\frak n}$ and $\frak n$ induce the same isomorphism
$\nu$ of $K$ into $\tilde{K^\prime}$.  By the first part of the proof
we have
$$\tilde{\frak m}_\mu\tilde{\frak n}=\frak D_{\mu, \nu}$$
for every isomorphism $\mu$ in $\underline\mu$, and hence
$$\frak m\frak n=\sum_{\underline\mu}\frak D_{\mu, \nu}.$$

  (3).  We finally consider the case when one of $\frak m$
and $\frak n$ is a $K$-unary prime  divisor.  Assume first that
$\frak m$ is a $K$-unary prime  divisor, then by Lemma  4.6
$\frak m\frak n=\frak m\nu$.  On the other hand, we have
$\mu\bar{\frak o}=\frak m$.  When $\frak n$ is $K$-unary,
we have $\nu\bar{\frak o}=\frak n$, and hence
$\frak D_{\mu, \nu}=0$.  When $\frak n$ is $K$-unary, we also have
$\frak m\frak n=0$ by the definition of divisor residues, and 
therefore $\frak m\frak n=\frak D_{\mu, \nu}$.
 When $\frak n$ is binary, we have
$\nu\bar{\frak o}=\Pi_{\Delta\frak n/K\nu}(\bar{\frak o})
\nu^{-1}$.  It follows from definition that
$\frak D_{\mu, \nu}=\frak m\nu$.  Therefore, we have
$\frak m\frak n=\frak D_{\mu, \nu}$.  

  Next, assume that $\frak n$ is $K$-unary and that $\frak m$
is binary.  Since $\frak m\frak n=\frak m(\frak n)$ by Theorem 4.11
(with the role of $K$ and $K^\prime$ being interchanged),
we have
$\frak m\frak n=\Pi_{\Delta\frak m/K^\prime}(\frak n\mu_0)$.
On the other hand, we have
$\frak D_{\mu, \nu}=\frak n\mu$ for every isomorphism
$\mu$ in $\underline\mu$.  Hence, we have
$$\sum_{\mu\in\underline\mu}\frak D_{\mu, \nu}=
\sum_{\mu\in\underline\mu}\frak n\mu=
\Pi_{\Delta\frak m/K^\prime}(\frak n\mu_0)=\frak m\frak n$$
by considering first the case when $N_{\Delta/K^\prime}(\frak m)=1$
and then the case when $N_{\Delta/K^\prime}(\frak m)$ is
arbitrary. \qed \enddemo

  In the following theorem, we take the norm of divisor residues.

\proclaim{Theorem 5.3}  Let $\frak m$ and $\frak n$ 
be two different prime
 divisors of $\Delta$, and let $(\mu, \mu^\prime)$ and $(\nu,
\nu^\prime)$ be the isomorphism pairs of $K^\prime$ induced 
by $\frak m$ and $\frak n$, respectively.
If $\frak m$ is a $K^\prime$-unary prime  divisor, define 
$\frak N_{\Delta\frak m/K^\prime\mu^\prime}(\frak n\frak m)
{\mu^\prime}^{-1}=N_{\Delta/K^\prime}(\frak n)\frak m$,
and if $\frak n$ is a $K^\prime$-unary prime  divisor, define 
$\frak N_{\Delta\frak n/K^\prime\nu^\prime}(\frak m\frak n)
{\nu^\prime}^{-1}$ $=N_{\Delta/K^\prime}(\frak m)\frak n$,
then
$$\frak N_{\Delta\frak m/K^\prime\mu^\prime}(\frak n\frak m)
{\mu^\prime}^{-1}=\frak N_{\Delta\frak n/K^\prime\nu^\prime}
(\frak m\frak n) {\nu^\prime}^{-1}.$$ 
If $\frak m$ is a $K$-unary prime  divisor, define 
$\frak N_{\Delta\frak m/K\mu}(\frak n\frak m)
{\mu}^{-1}=N_{\Delta/K}(\frak n)\frak m$,
and if $\frak n$ is a $K$-unary prime  divisor, define 
$\frak N_{\Delta\frak n/K\nu}(\frak m\frak n)
{\nu}^{-1}$ $=N_{\Delta/K}(\frak m)\frak n$,
then
$$\frak N_{\Delta\frak m/K\mu}(\frak n\frak m)
\mu^{-1}=\frak N_{\Delta\frak n/K\nu}(\frak m\frak n)\nu^{-1}.$$ 
\endproclaim

  \demo{Proof}  We first consider the case when $\frak m$ and
$\frak n$ are not $K^\prime$-unary.  Assume that
$\Delta\frak m=K\mu_0\cdot K^\prime$ and $\Delta 
\frak n=K\nu_0\cdot K^\prime$.  Let $\underline\mu$, $\underline\nu$
be isomorphism systems of $K$ coordinated to $\mu_0$, $\nu_0$.
By Lemma 5.2 we have $\frak m\frak n=\sum_{\mu\in\underline\mu}
\frak D_{\mu, \nu_0}$.  If $\sigma$ is an isomorphism of
$\Delta\frak n$ over $K^\prime$, then 
$(\frak m\frak n)\sigma=\sum_{\mu\in\underline\mu}
\frak D_{\mu, \nu_0\sigma}$.  It follows that
$$\frak N_{\Delta\frak n/K^\prime}(\frak m\frak n)=
\sum_{\mu\in\underline\mu, \nu\in\underline\nu}
\frak D_{\mu, \nu}.$$
Since $\frak D_{\mu, \nu}=\frak D_{\nu, \mu}$, we have
$$\frak N_{\Delta\frak n/K^\prime}(\frak m\frak n)=
\frak N_{\Delta\frak m/K^\prime}(\frak n\frak m).$$
This implies that
$$\frak N_{\Delta\frak m/K^\prime\mu^\prime}(\frak n\frak m)
{\mu^\prime}^{-1}=\frak N_{\Delta\frak n/K^\prime\nu^\prime}
(\frak m\frak n) {\nu^\prime}^{-1}$$
when $\Delta\frak n$ and $\Delta\frak m$
are not $K^\prime$-normalized.  
When $\frak m$ and $\frak n$ are not $K$-unary, 
a similar argument shows that the identity
$$\frak N_{\Delta\frak m/K\mu}(\frak n\frak m)
\mu^{-1}=\frak N_{\Delta\frak n/K\nu}(\frak m\frak n)\nu^{-1}$$
holds. 

  We now consider the case when one of $\frak m$ and $\frak n$
is $K^\prime$-unary, say $\frak m$.  Then
$\frak N_{\Delta\frak m/K^\prime\mu^\prime}(\frak n\frak m)
{\mu^\prime}^{-1}=N_{\Delta/K^\prime}(\frak n)\frak m$
by assumption.  If $\frak n$ is not $K^\prime$-unary, then
by Lemma  4.6 we have  $\frak m\frak n=\frak m\nu^\prime$.
It follows that 
$\frak N_{\Delta\frak n/K^\prime\nu^\prime}(\frak m\frak n)
{\nu^\prime}^{-1}=[\Delta\frak n:K^\prime\nu^\prime]
\frak m=N_{\Delta/K^\prime}(\frak n)\frak m$.
Therefore,
$$\frak N_{\Delta\frak m/K^\prime\mu^\prime}(\frak n\frak m)
{\mu^\prime}^{-1}=\frak N_{\Delta\frak n/K^\prime\nu^\prime}
(\frak m\frak n) {\nu^\prime}^{-1}$$
when $\frak n$ is not a $K^\prime$-unary prime  divisor.
If $\frak n$ is $K^\prime$-unary, then 
$N_{\Delta/K^\prime}(\frak n)=0$ and
$N_{\Delta/K^\prime}(\frak m)=0$ by definition.
Since 
$\frak N_{\Delta\frak n/K^\prime\nu^\prime}(\frak m\frak n)
{\nu^\prime}^{-1}$ $=N_{\Delta/K^\prime}(\frak m)\frak n$
by assumption, we have
$$\frak N_{\Delta\frak m/K^\prime\mu^\prime}(\frak n\frak m)
{\mu^\prime}^{-1}=\frak N_{\Delta\frak n/K^\prime\nu^\prime}
(\frak m\frak n) {\nu^\prime}^{-1}=0$$
when $\frak n$ is a $K^\prime$-unary prime  divisor.

  We finally consider the case when one of $\frak m$ and $\frak n$
is $K$-unary, say $\frak n$.  Then
$\frak N_{\Delta\frak n/K\nu}(\frak m\frak n)
{\nu}^{-1}$ $=N_{\Delta/K}(\frak m)\frak n$ by assumption.  
If $\frak m$ is not $K$-unary, then
by Lemma  4.6 we have  $\frak n\frak m=\frak n\mu$.
It follows that 
$$\frak N_{\Delta\frak m/K\mu}(\frak n\frak m)
\mu^{-1}=[\Delta\frak m:K\mu] \frak n=
N_{\Delta/K}(\frak m)\frak n.$$
Therefore,
$$\frak N_{\Delta\frak m/K\mu}(\frak n\frak m)
\mu^{-1}=\frak N_{\Delta\frak n/K\nu}
(\frak m\frak n) \nu^{-1}$$
when $\frak m$ is not a $K$-unary prime  divisor.
If $\frak m$ is $K$-unary, then 
$N_{\Delta/K}(\frak m)=0$ and
$N_{\Delta/K}(\frak n)=0$ by definition.
Since 
$\frak N_{\Delta\frak m/K\mu}(\frak n\frak m)
\mu^{-1}$ $=N_{\Delta/K}(\frak n)\frak m$
by assumption, we have
$$\frak N_{\Delta\frak m/K\mu}(\frak n\frak m)
\mu^{-1}=\frak N_{\Delta\frak n/K\nu}
(\frak m\frak n) \nu^{-1}=0$$
when $\frak n$ is a $K$-unary prime  divisor.

  This completes the proof of the theorem.  \enddemo

  Let $\frak A$ and $\frak b$ be divisors of $\Delta$, relatively
prime to each other.  Write $\frak A=\sum a_{\frak m}\frak m$ and
$\frak b=\sum b_{\frak n}\frak n$ as linear combination of prime
divisors of $\Delta$ with integer coefficients.  Define
$$\langle\frak A,\frak b\rangle_f=\sum_{\frak m, \frak n}a_{\frak m}
b_{\frak n}\frak N_{\Delta\frak n/K\nu}(\frak m\frak n)\nu^{-1}$$
and 
$$\langle\frak A,\frak b\rangle_f^\prime=\sum_{\frak m, \frak n}
a_{\frak m}b_{\frak n}\frak N_{\Delta\frak n/K^\prime\nu^\prime}
(\frak m\frak n){\nu^\prime}^{-1}.\tag 5.4$$
Note that $\langle\frak A,\frak b\rangle_f$ is a divisor of $K$ and that
$\langle\frak A,\frak b\rangle_f^\prime$ is a divisor of $K^\prime$.
It follows from Theorem 5.3 that $\langle\frak A,\frak b\rangle_f
=\langle\frak b, \frak A\rangle_f$ and 
$\langle\frak A,\frak b\rangle_f^\prime
=\langle\frak b, \frak A\rangle_f^\prime$.

  \proclaim{Lemma  5.5}  Let $\frak m$ be a binary prime  divisor of
$\Delta$ with $N_{\Delta/K^\prime}(\frak m)=1$.
Assume that $x$ is a separating element of $K$.  Denote by $dx$ the 
divisor $\frak D_x/u_x^2$, where $\frak D_x$ is the different 
of $K$ over $k(x)$ and $u_x$ is the denominator of $x$.  Then
$\langle\frak m, \frak m\rangle_f+ \langle dx, \frak m\rangle_f$
and $\langle\frak m, \frak m\rangle_f^\prime+ \langle dx, 
\frak m\rangle_f^\prime$ are principal.  \endproclaim

  \demo{Proof}  Denote by $\mu$, $\mu^\prime$ the isomorphism pair
of $K$, $K^\prime$ induced by $\frak m$.  Assume that 
$\Delta\frak m$ is $K^\prime$-normalized.  Since
$N_{\Delta/K^\prime}(\frak m)=1$, we have 
$K\mu\subseteq \Delta\frak m=K^\prime$.  
Since $x$ is an separating element of $K$, it is
integral for $\frak m$, and $x-x\mu\equiv 0$ modulo $\frak m$.
Denote by $\frak m_0$ the numerator of $x-x\mu$.  Then $\frak m_0$ is
a prime  divisor of $K^\prime(x)$ of $K^\prime$-degree one, 
which is not a prime
 divisor of $k(x)$.  Hence, $\frak m_0$ is not contained in 
the discriminant of $\Delta/K^\prime(x)$, which consists of only prime
 divisors of $k(x)$.  It follows that $\frak m_0$ is unramified 
in $\Delta$ over
$K^\prime(x)$.  This implies that $\frak m_0$ contains $\frak m$ only once.
It follows that there exists a  divisor $\frak n$, prime
to $\frak m$, such that $(x-x\mu)=\frak m-\frak n$.  By Theorem 4.8,
$\langle\frak m, \frak m\rangle_f- \langle \frak n, \frak m\rangle_f$
and $\langle\frak m, \frak m\rangle_f^\prime- \langle \frak n, 
\frak m\rangle_f^\prime$ are principal.
  Thus, the stated result is equivalent to the
statement that $\langle\frak m, \frak m\rangle_f
+ \langle dx, \frak m\rangle_f$
and $\langle\frak m, \frak m\rangle_f^\prime+ \langle dx, 
\frak m\rangle_f^\prime$ are principal.  
Since $\Delta\frak m=K^\prime$ and
since $dx$ is $K$-unary, by Lemma  4.6 we have
$$\aligned \langle\frak n, \frak m\rangle_f
+ \langle dx, \frak m\rangle_f
&=\frak N_{K^\prime/K\mu}(\frak n\frak m)\mu^{-1}+
\frak N_{K^\prime/K\mu}((dx)\mu)\mu^{-1}\\
&=\frak N_{K^\prime/K\mu}\left(\frak n\frak m+(dx)\mu
\right)\mu^{-1}\endaligned $$
and
$$\langle\frak n, \frak m\rangle_f^\prime+ \langle dx, 
\frak m\rangle_f^\prime=\frak n\frak m+(dx)\mu.$$
Therefore, it suffices to show that $\frak n\frak m+(dx)\mu$ is 
principal.

  Let $\tilde{K^\prime}$ be a splitting extension for $\frak m_0$.
Then $\frak m_0$ is unramified in $\tilde \Delta$
over $\tilde{K^\prime}(x)$, where $\tilde{\Delta}$ is the double
field of $K$, $\tilde{K^\prime}$.  Since $\Delta\frak m=K^\prime$, 
$\frak m$ remains a prime  divisor of $\tilde\Delta$.
Thus, $\frak m_0$ has in $\tilde\Delta$ the decomposition
$$\frak m_0=\frak m+\sum_{\sigma\in \underline\sigma, \, \sigma\neq 1}
\tilde{\frak m_{\sigma}}$$
into prime  divisors $\tilde{\frak m_\sigma}$ with 
$N_{\tilde{\Delta}/\tilde{K^\prime}}(\tilde{\frak m_\sigma})=1$,
where $\underline\sigma$ is the isomorphism system of $K\mu$ into
$\tilde{K^\prime}$.  It follows that
$$(x-x\mu)=\frak m_0-(u_x+u_x\mu)=
\frak m+\sum_{\sigma\in \underline\sigma, \, \sigma\neq 1}
\tilde{\frak m_{\sigma}}-(u_x+u_x\mu).$$
This implies that
$$-\frak n=\sum_{\sigma\in \underline\sigma, \, \sigma\neq 1}
\tilde{\frak m_{\sigma}}-(u_x+u_x\mu),$$
and hence
$$-\frak n\frak m=\sum_{\sigma\in \underline\sigma, \, \sigma\neq 1}
\tilde{\frak m_{\sigma}}\frak m-(u_x+u_x\mu)\frak m.$$
It follows from the part (1) in the proof of Lemma 5.2 that
$$\tilde{\frak m_{\sigma}}\frak m=\frak D_{\mu\sigma, \mu}.$$
Since $(u_x+u_x\mu)\frak m=2u_x\mu$, in order to prove that
$\frak n\frak m+(dx)\mu$ is principal it suffices to show that
$$\frak D_x\mu-\sum_{\sigma\in \underline\sigma, \, \sigma\neq 1}
\frak D_{\mu\sigma, \mu}$$
is principal.  But, by the approximation theorem on $K$ and
[2, Theorem 6 on page 92] we have 
$$\frak D_x\mu=\sum_{\sigma\in \underline\sigma, \, \sigma\neq 1}
\frak D_{\mu\sigma, \mu}.\qed$$  \enddemo

  By Lemma 5.5 and its proof, we obtain the following useful result.

 \proclaim{Corollary 5.6}  Let $\frak m$ be a binary prime  divisor of
$\Delta$ with $N_{\Delta/K^\prime}(\frak m)=1$.
Denote by $\mu$, $\mu^\prime=1$ the isomorphism pair of $K$,
$K^\prime$ induced by $\frak m$.  
Assume that $x$ is a separating element of $K$.  
Then $\langle\frak m, \frak m\rangle_f+N_{\Delta/K}(\frak m)(dx)$
and $\langle\frak m, \frak m\rangle_f^\prime+(dx)\mu$ 
are principal.  \endproclaim

\heading
6.  A residue scalar product of divisors
\endheading

   When $\infty$ is an infinite place of $k$, we view $K^\prime$ as an
algebraic function field of one variable over $\Bbb C$.  Denote 
by $K^\prime_\infty$ the set of all places of $K^\prime$.  
If $f$ is an element of $K^\prime$ and if $v$ is a place 
of $K^\prime$, then either $v$ is a pole of
$f$, in which case we say that $f$ takes the value $\infty$ at
$v$, or this is not the case, and then the value $f(v)$ taken
by $f$ at $v$ (which is the residue class of $f$ modulo $v$)
is a complex number.  Denote by $\hat\Bbb C$ the Riemann sphere,
obtained by adjunction of a point $\infty$ to $\Bbb C$.
Consider $K^\prime_\infty $ as a topological space whose topology
is the weakest topology with respect to which the mappings 
$f$: $v\to f(v)$
of $K^\prime_\infty$ into $\hat\Bbb C$ are continuous.  Then 
$K^\prime_\infty$ is a compact Riemann surface; see [4, Chapter VII].

   Let $\frak A$ be a divisor of $K^\prime$.  Then the extension 
of $\frak A$ to $K^\prime_\infty$ is a divisor of $K^\prime_\infty$, 
and there exists a line bundle $L$ with a Hermitian metric $\|\cdot\|$ 
and with a meromorphic section $s$ whose divisor is the extension 
of $\frak A$ to $K^\prime_\infty$ (cf. [1, \S2]).  
We assume that the metric $\|\cdot\|$ on $L$ satisfies
$$(\deg L)d\mu^\prime_\infty=c_1(\|\cdot\|, L),$$
where $c_1(\|\cdot\|, L)$ is the Chern form of the metric 
$\|\cdot\|$ on the line bundle $L$.  Define
$$\frak A_{K^\prime}=\frak A
+\sum_\infty \left(-\int_{K^\prime_\infty}\log \|s\| 
d\mu^\prime_\infty\right) K^\prime_\infty, $$
where the sum is over the infinite places of $k$.
Then $\frak A_{K^\prime}$ is an Arakelov divisor of $K^\prime$. 
We call $\frak A_{K^\prime}$ the Arakelov divisor
of $K^\prime$ obtained from the divisor $\frak A$. 
Let $\frak A$ and $\frak b$ be divisors of $\Delta$, relatively prime
to each other.  By (5.4), $\langle\frak A, \frak b\rangle_f^\prime$ 
is a divisor of $K^\prime$, and hence we can obtain an Arakelov divisor 
$\langle\frak A, \frak b\rangle_{K^\prime}$
of $K^\prime$ from the divisor $\langle\frak A, \frak b\rangle_f^\prime$.
By (5.4), we have $\langle\frak A, \frak b\rangle_f^\prime
=\langle\frak b, \frak A\rangle_f^\prime$.  This implies that 
$\langle\frak A, \frak b\rangle_{K^\prime}
=\langle\frak b, \frak A\rangle_{K^\prime}$.   If $f$ is an element 
of $K^\prime$, we put 
$v_\infty (f)=-\int_{K^\prime_\infty}\log |f|d\mu^\prime_\infty$.
Let $(f)_{\text{fin}}$ be the divisor of $f$ on $K^\prime$.  
Then a principal Arakelov divisor is given by
$$(f)_{K^\prime}=(f)_{\text{fin}}+\sum v_\infty(f)K^\prime_\infty$$
where the sum is over the infinite places of $k$.
It is clear that $\frak A_{K^\prime}$ is a principal Arakelov 
divisor of $K^\prime$ if $\frak A$ is a principal divisor of $K^\prime$.
Hence, by Theorem 4.8, $\langle\frak A, \frak b\rangle_{K^\prime}$ 
is a principal Arakelov divisor of $K^\prime$ if $\frak A$ is 
a principal divisor of $\Delta$.
 
 Let $x^\prime$ be a separating element of $K^\prime$.  Denote  
by $dx^\prime$ the divisor $\frak D_{x^\prime}/u_{x^\prime}^2$, where 
$\frak D_{x^\prime}$ is the different of $K^\prime$ over $k(x^\prime)$ 
and $u_{x^\prime}$ is the denominator of $x^\prime$.  Then $dx^\prime$
belongs to the class of divisors of differentials of $K^\prime$ (cf.
[5, Chapter 25]).   Let $\frak d_{K^\prime}$ be the Arakelov 
divisor of $K^\prime$ obtained from the divisor $dx^\prime$.   
For every Arakelov divisor $D$ of $K^\prime$, we define
$$\deg_{K^\prime} D=(D\cdot\frak d_{K^\prime}), \tag 6.1$$
where $(\,\,\cdot\,\,)$ is the Arakelov intersection product [1, \S1].  
Let $\frak A$ and $\frak b$ be divisors of $\Delta$, relatively prime 
to each other.  We define
$$\langle\frak A, \frak b\rangle=\deg_{K^\prime}\langle\frak A, 
\frak b\rangle_{K^\prime}. \tag 6.2$$

   If $\frak A$ is a divisor of $\Delta$, we denote by 
$\frak A|K^\prime$ 
the Arakelov divisor of $K^\prime$ obtained from the restriction of 
$\frak A$ to $K^\prime$.    If $\frak m, \frak n$ are two different 
prime  divisors of $\Delta$, we define
$$\{ \frak m, \frak n\}= N_{\Delta/K^\prime}(\frak n)
\deg_{K^\prime}(\frak m|K^\prime)
+N_{\Delta/K^\prime}(\frak m)\deg_{K^\prime}(\frak n|K^\prime).$$
Similarly, we define $\{\frak A, \frak b\}$ for all divisors 
$\frak A$, $\frak b$ of $\Delta$ by linearity.

 \proclaim{Definition 6.3}  Let $\frak A, \frak b$ be divisors 
of $\Delta$ which are relatively prime to each other.  
Then a residue scalar product 
$\langle\frak A, \frak b\rangle_r$ of $\frak A$ 
and $\frak b$ is defined by
$$\langle \frak A, \frak b\rangle_r
=\{ \frak A, \frak b\}-\langle \frak A, \frak b\rangle.$$
 \endproclaim

  If $(f)$ is a principal divisor of $\Delta$, 
it follows from (6.1) and Theorem 4.8 that 
$\langle (f), \frak A\rangle_r=0$ for any  divisor
$\frak A$ of $\Delta$, and hence the residue scalar product 
$\langle \frak A, \frak b\rangle_r$ 
is well-defined on the classes of divisors of 
$\Delta$ modulo principal divisors.  By (5.4) and Theorem 5.3,
the residue scalar product is bilinear and symmetric.  

  Part of Roquette's theory is to prove that the residue 
scalar product of a divisor of the double field with itself 
is nonnegative for function fields over a finite constants
field by using the Riemann-Roch theorem (cf. Li [7]).  
Then the Riemann hypothesis for function fields over a 
finite constants field follows from the Schwarz inequality.
Searching for the analogue for number fields of Roquette's
proof of the Riemann hypothesis for function fields over a
finite constants field, we were led to the sequence of numbers 
whose positivity is equivalent to the Riemann hypotheses [6].
We conjecture that the residue scalar product, which we 
constructed for function fields over a number field, of a 
divisor with itself is nonnegative.  This nonnegativity,
if true, may be related to that of the sequence of numbers 
found in Li [6].

\enddocument